\documentclass{smfart}
\usepackage{lmodern}
\usepackage[OT2,T1]{fontenc}
\usepackage[francais]{babel}
\usepackage{amsmath,amssymb,amsfonts,amsthm}
\usepackage[pdfborder={0 0 0}]{hyperref}

\advance\textheight -1.1mm
\advance\topmargin .55mm

\usepackage[ps,dvips,arrow,matrix,tips,line,curve]{xy}
\entrymodifiers={+!!<0pt,\the\fontdimen22\textfont2>}
\SelectTips{cm}{10}

\title[Groupe de Brauer et points entiers]{Groupe de Brauer et points entiers de deux familles de surfaces cubiques affines}

\date{30 juillet 2010}

\advance\abstractmargin.15cm
\begin{abstract}
Soit~$a$ un entier non nul.
Si~$a$ n'est pas de la forme $9n\pm 4$ pour un $n \in \Z$, il n'y
a pas d'obstruction de Brauer--Manin \`a l'existence d'une solution de
l'\'equation $x^3+y^3+z^3=a$ en entiers $x, y, z \in \Z$.
D'autre part, il n'y a pas d'obstruction de Brauer--Manin \`a l'existence d'une
solution de l'\'equation $x^3+y^3+2z^3=a$ en entiers $x, y, z \in \Z$.
\end{abstract}

\author[J.-L. Colliot-Th\'el\`ene]{Jean-Louis Colliot-Th\'el\`ene}
\email{jlct@math.u-psud.fr}
\address{C.N.R.S., U.M.R.~8628\\
Math\'ematiques, B\^atiment 425\\
Universit\'e Paris-Sud\\
F-91405 Orsay\\
France}

\author[O. Wittenberg]{Olivier Wittenberg}
\email{wittenberg@dma.ens.fr}
\address{D\'epartement de math\'ematiques et applications\\
\'Ecole normale sup\'erieure\\
45~rue d'Ulm\\
75230 Paris Cedex 05\\
France}

\newcommand{\Gm}{\mathbf{G}_\mathrm{m}}
\newcommand{\sU}{{\mathcal{U}}}
\newcommand{\sV}{{\mathcal{V}}}
\newcommand{\sM}{{\mathcal{M}}}
\newcommand{\sD}{{\mathcal{D}}}
\newcommand{\sH}{{\mathcal{H}}}
\newcommand{\sX}{{\mathcal{X}}}
\newcommand{\sZ}{{\mathcal{Z}}}
\newcommand{\kbar}{{\mkern1mu\overline{\mkern-1mu{}k\mkern-1mu}\mkern1mu}}
\newcommand{\Fbar}{{\mkern.5mu\overline{\mkern-.5mu{}\F\mkern-2mu}\mkern2mu}}
\newcommand{\mmu}{{\boldsymbol{\mu}}}
\newcommand{\tors}[2]{{\vphantom{#2}}_{#1}{#2}}
\newcommand{\Cbarre}{{\mkern2.5mu\overline{\mkern-2.5mu{}C\mkern-1mu}\mkern1mu}}
\newcommand{\Dbarre}{{\mkern2.5mu\overline{\mkern-2.5mu{}D\mkern-1mu}\mkern1mu}}
\newcommand{\Qbarre}{{\mkern1mu\overline{\mkern-1mu{}\Q\mkern-1mu}\mkern1mu}}
\newcommand{\Xbarre}{{\mkern2.5mu\overline{\mkern-2.5mu{}X\mkern-1mu}\mkern1mu}}
\newcommand{\Ubarre}{{\mkern2.5mu\overline{\mkern-2.5mu{}U\mkern-1mu}\mkern1mu}}
\newcommand{\RGamma}{\mathrm{R}\Gamma}
\newcommand{\Spec}{\mathop\mathrm{Spec}\nolimits}
\newcommand{\Gal}{\mathop\mathrm{Gal}\nolimits}
\newcommand{\Pic}{\mathop\mathrm{Pic}\nolimits}
\newcommand{\Br}{\mathop\mathrm{Br}\nolimits}
\newcommand{\Hom}{\mathop\mathrm{Hom}\nolimits}
\newcommand{\Cores}{\mathop\mathrm{Cores}\nolimits}
\newcommand{\Ker}{\mathop\mathrm{Ker}\nolimits}
\newcommand{\Res}{{\mathrm{Res}}}
\newcommand{\A}{{\mathbf A}}
\newcommand{\Q}{{\mathbf Q}}

\newcommand{\F}{{\mathbf F}}
\newcommand{\Fp}{\F_{\mkern-2mup}}
\newcommand{\Fbarp}{\Fbar_{\mkern-2mup}}
\newcommand{\Z}{{\mathbf Z}}
\newcommand{\N}{{\mathbf N}}
\newcommand{\R}{{\mathbf R}}
\newcommand{\congru}[3]{#1 \equiv #2 \!\!\mod #3}
\newcommand{\pascongru}[3]{#1 \not\equiv #2 \!\!\mod #3}
\newcommand{\nr}{{\mathrm{nr}}}
\newcommand{\uplet}[2]{#1, \mskip2.5mu \ldots \mskip-1mu, \mskip2.5mu #2}
\newcommand{\ev}{{\mathrm{ev}}}
\renewcommand{\P}{{\mathbf P}}
\renewcommand{\leq}{\leqslant}
\renewcommand{\geq}{\geqslant}
\renewcommand{\emptyset}{\varnothing}
\DeclareMathOperator{\inv}{inv}
\newcommand{\isoto}{\myxrightarrow{\,\sim\,}}
\makeatletter
\def\myrightarrow{{\setbox\z@\hbox{$\rightarrow$}\dimen0\ht\z@\multiply\dimen0 6\divide\dimen0 10\ht\z@\dimen0\box\z@}}
\def\myrightarrowfill@{\arrowfill@\relbar\relbar\myrightarrow}
\newcommand{\myxrightarrow}[2][]{\ext@arrow 0359\myrightarrowfill@{#1}{#2}}
\makeatother

\numberwithin{equation}{section}
\theoremstyle{plain}
\newtheorem{prop}{Proposition}[section]
\newtheorem{thm}[prop]{Th\'eor\`eme}
\newtheorem{lem}[prop]{Lemme}
\theoremstyle{remark}
\newtheorem{rmq}[prop]{Remarque}
\newtheorem{rmqs}[prop]{Remarques}
\newtheorem{questions}[prop]{Questions}
\newtheorem{exemple}[prop]{Exemple}
{\setbox0\hbox{$ $}}\fontdimen16\textfont2=\fontdimen17\textfont2

\makeatletter
\def\@setaddresses{\par\nobreak
  \begingroup
  \parindent-2em\leftskip2em
  \rightskip=0pt plus 20pt
  \emergencystretch .5\textwidth
  \exhyphenpenalty=-100
  \interlinepenalty\@M
  \def\baselinestretch{1}\normalfont\footnotesize
  \def\\{\unskip, \penalty-10\ignorespaces}%
  \def\cond@bullet {\unskip
         {\discretionary{}{}{\hbox{\ $\bullet$\ }}}}
  \def\author##1{\ifhmode\par\nobreak \vskip\smallskipamount\fi
      {\scshape ##1}\let\address\firstaddress}%
  \def\firstaddress##1##2{\unskip, \let\address\otheraddress
         \penalty-20\ignorespaces##2}%
  \def\otheraddress##1##2{\cond@bullet \ignorespaces##2}%
  \def\curraddr{\address}
  \let\address\firstaddress
  \def\email##1##2{\@ifnotempty{##2}%
        {\cond@bullet
         \hbox{\itshape Courriel~:}~{\ttfamily\ignorespaces ##2}}}%
  \def\urladdr##1##2{\@ifnotempty{##2}%
        {\cond@bullet
         {\itshape Url~:}~{\ttfamily\ignorespaces ##2}\par}}%
  \addresses
  \par\endgroup
}
\makeatother

\begin{document}
\maketitle 

\section{Introduction}

Il est connu depuis Ryley~\cite{ryley} que tout entier, et m\^eme tout nombre
rationnel, peut s'\'ecrire comme somme de trois cubes de nombres rationnels.
La question de savoir quels entiers s'\'ecrivent comme sommes de trois cubes
d'entiers relatifs, en revanche, est toujours ouverte.
Un tel entier ne peut \^etre congru \`a~$4$ ou~$5$ modulo~$9$.  Plusieurs auteurs
ont conjectur\'e que r\'eciproquement, tout entier non congru \`a~$4$ ou~$5$
modulo~$9$ est la somme de trois cubes d'entiers relatifs
(cf.~\cite[p.~623]{HBdensity}, \cite{connvaserstein}).

Un probl\`eme apparent\'e, remontant semble-t-il \`a Mordell (cf.~\cite{ko}),
consiste \`a d\'eterminer quels entiers s'\'ecrivent sous la forme $x^3+y^3+2z^3$
avec $x,y,z \in \Z$.  Il~se pourrait que tout entier, sans exception, poss\`ede
cette propri\'et\'e.  La relation bien connue
\begin{align}
\label{eq:identite-6a}
\left(a + 1\right)^3 + \left(a - 1\right)^3 - 2a^3 = 6a
\end{align}
montre au moins que les entiers multiples de~$6$ admettent une telle \'ecriture.

Il n'est pas exclu que des identit\'es similaires \`a~(\ref{eq:identite-6a})
puissent suffire \`a d\'emontrer que tout entier non congru \`a~$4$ ou~$5$
modulo~$9$ est une somme de trois cubes.
Toutefois, Vaserstein~\cite{vaserstein} a \'etabli que si
$x(t),y(t),z(t) \in \Q[t]$
sont des polyn\^omes tels que $x(t)^3+y(t)^3+z(t)^3=t$,
alors l'un de $x(t)$, $y(t)$ et $z(t)$ est n\'ecessairement de degr\'e~$\geq 5$.

La question de la r\'esolubilit\'e des \'equations $x^3+y^3+z^3=a$ et
$x^3+y^3+2z^3=a$ en nombres entiers relatifs, pour~$a$ fix\'e, est signal\'ee dans
le recueil de Guy~\cite[D5, p.~231]{guy} et a donn\'e lieu \`a de nombreuses
exp\'eriences sur ordinateur, qui jusqu'ici n'ont permis d'exhiber aucun triplet
$(x,y,z)\in \Z^3$ tel que $x^3+y^3+z^3=33$ ou $x^3+y^3+2z^3=148$
(cf.~\cite{elsenhansjahnel}, \cite{koyama}).

\'Etant donn\'e un sch\'ema~$\sU$ de type fini sur~$\Z$, il r\'esulte de la th\'eorie du
corps de classes global que l'image de l'application diagonale
$\sU(\Z) \to \prod \sU(\Z_p)$, o\`u~$p$ parcourt l'ensemble des places de~$\Q$
(en convenant que $\Z_\infty=\R$), est incluse dans le sous-ensemble des
points ad\'eliques orthogonaux au groupe de Brauer de $U = \sU \otimes_\Z \Q$
pour l'accouplement de Brauer--Manin (cf.~\cite[\textsection1]{ctxu}).
La vacuit\'e de ce sous-ensemble permet dans certains cas d'expliquer celle
de~$\sU(\Z)$ (cf.~\cite{ctxu}, \cite{KT}).  On parle alors d'obstruction de
Brauer--Manin enti\`ere (ou simplement d'obstruction de Brauer--Manin) \`a l'existence d'un $\Z$\nobreakdash-point de~$\sU$.
C'est un cas particulier de l'obstruction de Brauer--Manin \`a
l'approximation forte sur~$U$.

Cassels~\cite{casselsmathcomp} s'\'etait rendu compte que la loi de r\'eciprocit\'e
cubique impose une contrainte non triviale aux solutions enti\`eres de
l'\'equation $x^3+y^3+z^3=3$ (\`a savoir, les entiers~$x$, $y$ et~$z$ sont
n\'ecessairement congrus entre eux modulo~$9$).  Son argument revient en r\'ealit\'e
\`a exhiber une obstruction de Brauer--Manin \`a l'approximation forte sur la
surface affine d'\'equation $x^3+y^3+z^3=3$ sur~$\Q$. (Il s'agit m\^eme d'une
obstruction \`a l'approximation faible; voir la remarque~\ref{rq:cassels}
ci-dessous.)
Il est naturel de se demander si de fa\c{c}on plus g\'en\'erale, il se pourrait qu'une
obstruction de Brauer--Manin interdise, pour certains entiers~$a$,
l'existence de $x,y,z \in \Z$ tels que $x^3+y^3+z^3=a$ ou $x^3+y^3+2z^3=a$
(sous l'hypoth\`ese, pour la premi\`ere \'equation, que~$a$ n'est pas congru \`a~$4$
ou~$5$ modulo~$9$).  Le but de cet article est de d\'emontrer qu'il n'en est
rien: les diverses lois de r\'eciprocit\'e englob\'ees dans l'obstruction de
Brauer--Manin enti\`ere ne permettent d'exclure aucun entier~$a$.

Le texte est organis\'e comme suit.  Les paragraphes~\ref{par:corps-arbitraire}
et~\ref{par:rationnels} sont consacr\'es \`a la d\'etermination des groupes de
Brauer des surfaces affines d'\'equations $x^3+y^3+z^3=a$ et $x^3+y^3+2z^3=a$
sur~$\Q$ (propositions~\ref{prop:proj-br}, \ref{prop:brauer-affine}
et~\ref{prop:brauer-affine2}).  Les d\'emonstrations des
propositions~\ref{prop:brauer-affine} et~\ref{prop:brauer-affine2} requi\`erent
l'utilisation d'information arithm\'etique non triviale (classification des
courbes elliptiques sur~$\Q$ isog\`enes \`a une courbe elliptique donn\'ee) et ne
s'\'etendraient pas telles quelles \`a des \'equations plus g\'en\'erales ou \`a des corps
de nombres autres que~$\Q$.  Au paragraphe~\ref{par:sur-les-entiers}, tirant
parti de la structure de groupe des courbes elliptiques $x^3+y^3+z^3=0$ et
$x^3+y^3+2z^3=0$ sur~$\Fp$, nous \'etablissons le th\'eor\`eme principal de
l'article (th\'eor\`eme~\ref{thm:sur-les-entiers}).  Finalement, le
paragraphe~\ref{par:remarques-questions} rassemble divers compl\'ements \`a la
d\'emonstration du th\'eor\`eme~\ref{thm:sur-les-entiers}, dont certains pourraient
avoir un int\'er\^et ind\'ependant, discute quelques exemples
et soul\`eve la question de la densit\'e de~$U(k)$
dans l'ensemble des points ad\'eliques de~$U$ orthogonaux \`a~$\Br(U)$ (pour la
topologie ad\'elique modifi\'ee aux places infinies; voir la question~\ref{qu}~(ii)
pour un \'enonc\'e pr\'ecis) lorsque~$U$ est le compl\'ementaire d'une section
hyperplane lisse dans une surface cubique lisse sur un corps de nombres.

\bigskip
\emph{Remerciements.}
Des calculs partiels sur les surfaces $x^3+y^3+z^3=at^3$ avaient \'et\'e
r\'ealis\'es par Venapally Suresh il y a quelques ann\'ees.  Nous
savons gr\'e \`a Samir Siksek d'avoir attir\'e notre attention sur l'\'equation
$x^3+y^3+2z^3=a$.  Le premier auteur remercie l'Institut Alfr\'ed R\'enyi
(Budapest) et l'Institut Hausdorff (Bonn) pour leur hospitalit\'e.

\bigskip
\emph{Notation.}
Si~$k$ est un corps et $a$ un \'el\'ement de~$k$, on notera $k(\sqrt[3]a)$ une
extension minimale de~$k$ dans laquelle~$a$ soit un cube; le
symbole~$\sqrt[3]a$ d\'esignera une racine cubique de~$a$, fix\'ee une fois pour
toutes, dans cette extension.  Selon cette convention, l'extension
$k(\sqrt[3]a)/k$ est toujours de degr\'e~$1$ ou~$3$.

\section{Sur un corps arbitraire}
\label{par:corps-arbitraire}

Dans ce paragraphe nous d\'eterminons le groupe de Brauer des surfaces cubiques
projectives d'\'equations $x^3+y^3+z^3=at^3$ et $x^3+y^3+2z^3=at^3$ et en
exhibons des g\'en\'erateurs explicites (proposition~\ref{prop:proj-br}).
Le~calcul, purement alg\'ebrique, vaut sur un corps arbitraire de
caract\'eristique diff\'erente de~$3$ ne contenant pas de racine cubique primitive
de l'unit\'e.

Soient~$k$ un tel corps et~$X$ une vari\'et\'e g\'eom\'etriquement irr\'eductible et
lisse sur~$k$.  Notons $K=k[j]/(j^2+j+1)$ et $G=\Gal(K/k)=\{1,\sigma\}$,
o\`u~$\sigma$ est d\'efini par $\sigma(j)=j^2$.  Notons~$k(X)$ et~$K(X)$ les corps
de fonctions rationnelles de~$X$ et de $X \otimes_k K$.

Soit~$\mmu_3$ le $G$\nobreakdash-module $\{1,j,j^2\}$.  Pour $f,g\in K(X)^*$,
le cup-produit $\{f,g\}$ des classes de~$f$ et de~$g$ dans
$K(X)^*/K(X)^{*3}=H^1(K(X),\mmu_3)$ est un \'el\'ement de
$H^2(K(X),\mmu_3^{\otimes 2})$.  Nous noterons $(f,g)_j \in H^2(K(X),\mmu_3)$
l'image de $\{f,g\} \otimes j$ par l'isomorphisme $G$\nobreakdash-\'equivariant
canonique
$$
H^2(K(X),\mmu_3^{\otimes 2}) \otimes \mmu_3 \isoto H^2(K(X),\mmu_3^{\otimes 3}) = H^2(K(X),\mmu_3)
$$
(en remarquant que  $\mmu_3^{\otimes 3}=\mmu_3$ puisque
$\mmu_3^{\otimes 2}=\Z/3\Z$); ce symbole est $\Z$\nobreakdash-bilin\'eaire et
antisym\'etrique en~$f$ et~$g$.  Il r\'esulte des \'egalit\'es
$$
\sigma\mkern-1.5mu\left( \{a,b\}\otimes j \right) = \{\sigma(a),\sigma(b)\} \otimes \sigma(j) = \{\sigma(a),\sigma(b)\} \otimes j^2 = -\{\sigma(a),\sigma(b)\} \otimes j
$$
que l'action de~$G$ sur~$(a,b)_j$ est donn\'ee par
\begin{align}
\label{eq:action}
\sigma\mkern-2mu\left( (a,b)_j \right) = -(\sigma(a),\sigma(b))_j\text{.}
\end{align}

Rappelons enfin que le groupe $H^2(k(X),\mmu_3)$ (resp.~$H^2(K(X),\mmu_3)$)
s'identifie au sous-groupe de $3$\nobreakdash-torsion de~$\Br(k(X))$ (resp.~de
$\Br(K(X))$) et que la lissit\'e de~$X$ entra\^ine que les fl\`eches naturelles
$\Br(X) \to \Br(k(X))$ et $\Br(X\otimes_k K) \to \Br(K(X))$ sont des
injections.

\begin{prop}
\label{prop:proj-br}
Soit~$k$ un corps de caract\'eristique diff\'erente de~$3$ ne contenant pas de
racine cubique primitive de l'unit\'e.  Soit $a \in k^*$. Notons
$X \subset \P^3_k$ et $X' \subset \P^3_k$ les surfaces projectives et lisses
d'\'equations homog\`enes respectives $x^3+y^3+z^3=at^3$ et $x^3+y^3+2z^3=at^3$.
Les classes
\begin{align*}
A&=\Cores_{K(X)/k(X)} ((x+jy)/(x+y),a)_j \in H^2(k(X),\mmu_3) \subset \Br(k(X)) \rlap{\text{,}}\\
A'&=\Cores_{K(X')/k(X')} ((x+jy)/(x+y),4a)_j \in H^2(k(X'),\mmu_3) \subset \Br(k(X'))
\end{align*}
appartiennent aux sous-groupes $\Br(X) \subset \Br(k(X))$ et
$\Br(X') \subset \Br(k(X'))$.  Si~$a$ n'est pas un cube dans~$k$, le quotient
$\Br(X)/\Br(k)$ est d'ordre~$3$, engendr\'e par l'image de~$A$.  Si aucun
de~$a$, $2a$, $4a$ n'est un cube dans~$k$, le quotient $\Br(X')/\Br(k)$ est
d'ordre~$3$, engendr\'e par l'image de~$A'$.
\end{prop}

\begin{proof}
On suppose que~$a$ n'est pas un cube.  Il est connu
(cf.~\cite[\textsection45.3]{Manin}, \cite[p.~430]{CTS}) que le quotient
$\Br(X\otimes_k K)/\Br(K)$ est isomorphe \`a $(\Z/3)^2$, que les deux
\'el\'ements $$\alpha=\left((x+jy)/(x+y),a\right)_j$$
et $$\beta=\left((x+z)/(x+y),a\right)_j$$ de $\Br(K(X))$ appartiennent au
sous-groupe $\Br(X \otimes_k K)$, et que leurs images dans
$\Br(X\otimes_k K)/\Br(K)$ constituent une base de ce
$\Z/3\Z$\nobreakdash-espace vectoriel.  Ceci implique d\'ej\`a que
$A=\Cores_{K(X)/k(X)}(\alpha)$ appartient \`a~$\Br(X)$.

La formule~(\ref{eq:action}) montre que $\sigma(\beta)=-\beta$ et que
\begin{align}
\label{calcul:alpha}
\alpha-\sigma(\alpha)&=((x+jy)(x+j^2y)/(x+y)^2,a)_j=((x^3+y^3)/(x+y)^3,a)_j\\
\notag&=((x^3+y^3)/t^3,a)_j=(a-(z/t)^3,a)_j=((z/t)^3-a,a)_j=0
\end{align}
(cf.~\cite[Ch.~XIV, \textsection2, prop.~4]{corpslocaux}), d'o\`u
$\sigma(\alpha)=\alpha$.
Ainsi le groupe des invariants $\left(\Br(X \otimes_k K)/\Br(K)\right)^G$
est-il isomorphe \`a~$\Z/3\Z$, engendr\'e par~$\alpha$.

Les applications
\begin{align*}
\Res : \Br(X)/\Br(k) \longrightarrow \Br(X\otimes_k K)/\Br(K)\rlap{\text{,}} \\ \Cores :  \Br(X\otimes_k K)/\Br(K) \longrightarrow \Br(X)/\Br(k)
\end{align*}
de restriction et de corestriction satisfont $\Cores \circ \Res = 2$ et
$\Res \circ \Cores = 1 +\sigma$.
Comme la surface~$X\otimes_kk(\sqrt[3]{a}\mkern2mu)$
est $k(\sqrt[3]{a}\mkern2mu)$\nobreakdash-rationnelle
(cf.~\cite[\textsection13.7]{hardywright}),
le quotient
$\Br(X\otimes_kk(\sqrt[3]{a}\mkern2mu))/\Br(k(\sqrt[3]{a}\mkern2mu))$
est nul.  L'extension
$k(\sqrt[3]{a}\mkern2mu)/k$ \'etant en outre de degr\'e impair,
il s'ensuit que $\Br(X)/\Br(k)$ est d'ordre impair.
L'application~$\Res$ induit donc un isomorphisme
$$
\Br(X)/\Br(k) \isoto \left(\Br(X\otimes_kK)/\Br(K)\right)^G
$$
dont l'inverse est $-\Cores$. Compte tenu que la compos\'ee
$\Res \circ \Cores$ co\"incide sur
$\left(\Br(X \otimes_k K)/\Br(K)\right)^G$
avec la multiplication par~$2$, on conclut que $\Br(X)/\Br(k)$ est d'ordre~$3$
et est engendr\'e par $\Cores(\alpha)$.

Supposons maintenant que ni~$a$, ni~$2a$, ni~$4a$ ne soient des cubes
dans~$k$.  Si~$2$ est un cube dans~$k$, la surface~$X'$ est isomorphe \`a~$X$
et~$A'$ s'identifie \`a~$A$, de sorte qu'il n'y a rien de plus \`a d\'emontrer.
Supposons donc que~$2$ ne soit pas un cube.
D'apr\`es~\cite[\textsection1, proposition~1]{CTKS},
le quotient $\Br(X' \otimes_k K)/\Br(K)$ est alors isomorphe \`a~$\Z/3\Z$.
Comme d'autre part, notant~$\ell$ l'extension de degr\'e~$9$ de~$k$
donn\'ee par $\ell=k(\sqrt[3]{a},\sqrt[3]{2}\mkern2mu)$, la surface
$X'\otimes_k\ell$ est $\ell$\nobreakdash-rationnelle,
le groupe $\Br(X')/\Br(k)$ est d'ordre impair et l'application de restriction
induit \`a nouveau un isomorphisme
$$
\Br(X')/\Br(k) \isoto \left(\Br(X' \otimes_kK)/\Br(K)\right)^G
$$
dont l'inverse est l'oppos\'e de la corestriction.

L'\'el\'ement $\alpha'=\left((x+jy)/(x+y),4a\right)_j$ de $\Br(K(X'))$ appartient
au sous-groupe $\Br(X' \otimes_k K) \subset \Br(K(X'))$ car le diviseur de la
fonction rationnelle $(x+jy)/(x+y)$ sur $X' \otimes_k K$ est la norme d'un
diviseur sur $X' \otimes_k K(\sqrt[3]{4a}\mkern2mu)$
(cf.~\cite[Ex.~2.7.8~(c)]{CTS}).  Par cons\'equent
$A'=\Cores_{K(X')/k(X')}(\alpha')$ appartient \`a $\Br(X')$.

Pour conclure il suffit de v\'erifier que l'image de $\alpha'$ dans
$\Br(X' \otimes_kK)/\Br(K)$ est non nulle et est invariante par~$G$.
Qu'elle soit invariante par~$G$ r\'esulte d'un calcul similaire
\`a~(\ref{calcul:alpha}).  Elle est non nulle parce que son image dans
$\Br(X' \otimes_k K(\sqrt[3]2))/\Br(K(\sqrt[3]2))$ est non nulle: en effet,
sur le corps $K(\sqrt[3]2)$, les surfaces~$X$ et~$X'$ deviennent isomorphes
et~$\alpha'$ s'identifie \`a~$\alpha$, dont on a d\'ej\`a vu que l'image dans
$\Br(X \otimes_k K(\sqrt[3]2))/\Br(K(\sqrt[3]2))$ est non nulle.
\end{proof}

\section{Sur les rationnels}
\label{par:rationnels}

Le but de ce paragraphe est de d\'eterminer les groupes de Brauer des surfaces
cubiques affines d'\'equations $x^3+y^3+z^3=a$ et $x^3+y^3+2z^3=a$ sur le
corps~$\Q$ des rationnels (propositions~\ref{prop:brauer-affine}
et~\ref{prop:brauer-affine2}).

\begin{prop}
\label{prop:brauer-affine}
Soit $a \in \Q^*$.
Notons $X \subset \P^3_\Q$ la surface projective et lisse d'\'equation homog\`ene
$x^3+y^3+z^3=at^3$ et $U \subset \A^3_\Q$
la surface cubique affine d'\'equation $x^3+y^3+z^3=a$, de sorte que~$U$ est le
compl\'ementaire dans~$X$ de la section hyperplane d'\'equation $t=0$.  La fl\`eche
de restriction $\Br(X)\to\Br(U)$ est un isomorphisme.
\end{prop}

Sous les hypoth\`eses de la proposition~\ref{prop:brauer-affine}, si~$a$ n'est
pas un cube dans~$\Q$, le quotient $\Br(U)/\Br(\Q)$ est donc d'ordre~$3$,
engendr\'e par l'image de la classe~$A$ d\'efinie dans la
proposition~\ref{prop:proj-br}.

\begin{proof}
Notons $D \subset X$ la section hyperplane d'\'equation $t=0$.  Comme~$D$ est
lisse, il r\'esulte du th\'eor\`eme de puret\'e pour le groupe de Brauer que la suite
exacte de localisation s'\'ecrit
\begin{align}
\xymatrix{
0 \ar[r] & \Br(X) \ar[r] & \Br(U) \ar[r]^(.4){\partial_D} & H^1(D,\Q/\Z)
}
\end{align}
(cf.~\cite[\textsection6]{brauerIII}).  Pour \'etablir la proposition il suffit
donc de montrer que l'application $\partial_D$ est nulle.
Soient $A \in \Br(U)$ et $m=\partial_D(A) \in H^1(D,\Q/\Z)$.
Notons $P_0$, $P_1$, $P_2 \in X(\Q)$ les points de coordonn\'ees homog\`enes
$P_0=[1:-1:0:0]$, $P_1=[-1:0:1:0]$, $P_2=[0:-1:1:0]$.  Soit
$k=\Q(\sqrt[3]a\mkern2mu)$.  Soient $L_0$, $L_1$, $L_2 \subset \P^3_k$ les
droites d'\'equations
\begin{align*}
L_0 \colon x+y=0\text{, }z=\sqrt[3]a\mkern1mu{}t\rlap{\text{,}}\\
L_1 \colon x+z=0\text{, }y=\sqrt[3]a\mkern1mu{}t\rlap{\text{,}}\\
L_2 \colon y+z=0\text{, }x=\sqrt[3]a\mkern1mu{}t\rlap{\text{.}}
\end{align*}
Pour tout~$i$, la droite $L_i$ est contenue dans $X \otimes_\Q k$ et rencontre
$D \otimes_\Q k$ en un unique point, le point $P_i \otimes_\Q k$.
Notons $A_i \in \Br(L_i \cap U)$ la restriction de~$A$ \`a $L_i \cap U$.
Comme $L_i \cap U$ est isomorphe \`a $\A^1_k$ et que $\Br(\A^1_k)=\Br(k)$ (la
caract\'eristique de~$k$ \'etant nulle), la classe~$A_i$ est constante; son r\'esidu
au point $P_i \otimes_\Q k$ est donc nul.  Comme d'autre part $D\otimes_\Q k$
et~$L_i$ se rencontrent transversalement en $P_i\otimes_\Q k$, le r\'esidu
de~$A_i$ en $P_i\otimes_\Q k$ n'est autre que l'image de~$m$ par l'application
$H^1(D,\Q/\Z) \to H^1(k,\Q/\Z)$ d'\'evaluation en $P_i \otimes_\Q k$.  Il
s'ensuit que la classe $m(P_i) \in H^1(\Q,\Q/\Z)$ appartient pour tout~$i$ au
noyau de la fl\`eche de restriction $H^1(\Q,\Q/\Z) \to H^1(k,\Q/\Z)$.  Or cette
fl\`eche est injective puisque l'extension $k/\Q$ est soit triviale, soit
cubique et non cyclique.  D'o\`u la nullit\'e de $m(P_i)$ pour tout~$i$.
L'assertion~(b) du lemme ci-dessous permet d'en d\'eduire que $m=0$, ce qui
conclut la d\'emonstration de la proposition~\ref{prop:brauer-affine}.
\end{proof}

\begin{lem}
\label{lemme:cremona}
Soient $D \subset \P^2_\Q$ la courbe plane d'\'equation $x^3+y^3+z^3=0$
et $P_0$, $P_1$, $P_2 \in D(\Q)$ les points de coordonn\'ees
$P_0=[1:-1:0]$, $P_1=[-1:0:1]$, $P_2=[0:-1:1]$.  Alors:

(a) Le noyau de la fl\`eche $H^1(D,\Q/\Z)\to H^1(\Q,\Q/\Z)$ d'\'evaluation au
point~$P_0$ est isomorphe \`a $\Z/3\Z$.

(b) Un \'el\'ement de $H^1(D,\Q/\Z)$ nul en~$P_0$ et en~$P_1$ est
nul.
\end{lem}

\begin{proof}
Munissons~$D$ de la structure de courbe elliptique sur~$\Q$ ayant pour \'el\'ement
neutre le point rationnel~$P_0$.  Pour \'etablir le lemme~\ref{lemme:cremona} il
suffit de conna\^itre toute la torsion rationnelle pouvant appara\^itre dans une
courbe elliptique isog\`ene sur~$\Q$ \`a~$D$.  En effet, un \'el\'ement~$m$ de
$H^1(D,\Q/\Z)$ nul au point~$P_0$ s'identifie \`a une suite exacte
\begin{align}
\label{se:isog-cyclique}
\xymatrix{
0 \ar[r] & \Z/n\Z \ar[r] & E \ar[r] & D \ar[r] & 0
}
\end{align}
o\`u~$E$ est une courbe elliptique sur~$\Q$ et~$n$ est l'ordre de~$m$.
Pour $P \in D(\Q)$, la classe $m(P)$ est alors \'egale \`a l'image de~$P$ par la
fl\`eche $D(\Q) \to H^1(\Q,\Z/n\Z) \subset H^1(\Q,\Q/\Z)$ bord
de~(\ref{se:isog-cyclique}).

Le changement de variables $x=-4-3\eta$, $y=-4+3\eta$, $z=6\xi$, indiqu\'e par
Selmer dans~\cite{Selmer}, permet de ramener l'\'equation de~$D$ sous la forme
de Weierstrass $\eta^2=\xi^3-16/27$.  Cette courbe est not\'ee 27A1 dans les
tables de Cremona~\cite{cremona}.  Dans celles-ci, on lit de plus que la
classe d'isog\'enie de~$D$ comprend trois autres courbes elliptiques sur~$\Q$ et
que le sous-groupe de torsion du groupe de Mordell--Weil de chaque courbe
elliptique sur~$\Q$ isog\`ene \`a~$D$ est soit trivial soit d'ordre~$3$.  Par
cons\'equent une suite exacte de la forme~(\ref{se:isog-cyclique}) avec $n>1$ ne
peut exister que pour $n=3$.  Par dualit\'e, elle \'equivaut alors \`a la donn\'ee
d'une injection $\mmu_3 \hookrightarrow D$.

Les points de~$D$ annul\'es par~$3$ sont ceux qui v\'erifient l'\'equation $xyz=0$.
Comme l'intersection de~$D$ avec la droite $z=0$ est un sous-groupe
canoniquement isomorphe \`a~$\mmu_3$ et comme le sous-groupe $\{P_0,P_1,P_2\}$
est quant \`a lui canoniquement isomorphe \`a~$\Z/3\Z$, le sous-groupe de
$3$\nobreakdash-torsion de~$D$ se d\'ecompose en $\tors3D=\mmu_3 \oplus \Z/3\Z$.
Il~s'ensuit que $\Hom(\mmu_3,D)=\Hom(\mmu_3,\mmu_3)=\Z/3\Z$ (o\`u~$\Hom$ d\'esigne
l'ensemble des homomorphismes de groupes alg\'ebriques sur~$\Q$).

Le noyau de la fl\`eche $H^1(D,\Q/\Z)\to H^1(\Q,\Q/\Z)$ d'\'evaluation en~$P_0$
est donc isomorphe \`a $\Z/3\Z$, un g\'en\'erateur \'etant donn\'e par la suite exacte
$$
\xymatrix{
0 \ar[r] & \Z/3\Z \ar[r] & D/\mmu_3 \ar[r] & D \ar[r] & 0
}
$$
obtenue par passage au quotient \`a partir de la multiplication par~$3$ sur~$D$.
Si ce g\'en\'erateur s'annulait en~$P_1$, autrement dit s'il existait un point
rationnel de $D/\mmu_3$ s'envoyant sur~$P_1$, le groupe de Mordell--Weil de
$D/\mmu_3$ contiendrait un sous-groupe d'ordre~$9$.  Or tel n'est pas le cas
(par exemple parce que, comme on l'a vu, le groupe de Mordell--Weil d'une
courbe elliptique sur~$\Q$ isog\`ene \`a~$D$ ne contient jamais de sous-groupe
d'ordre~$9$), d'o\`u le lemme.
\end{proof}
 
\begin{rmqs}
\label{rq:HSC}
(i) De fa\c{c}on g\'en\'erale, si~$D$ est une vari\'et\'e propre, lisse et g\'eom\'etriquement
connexe sur un corps~$k$, si $\Dbarre$ est la vari\'et\'e obtenue par extension
des scalaires \`a une cl\^oture s\'eparable de~$k$ et si~$G$ d\'esigne le groupe de
Galois absolu de~$k$, la suite spectrale de Hochschild--Serre fournit une
suite exacte
$$
\xymatrix@C=1em{
0 \ar[r] & H^1(k,\Q/\Z) \ar[r] & H^1(D,\Q/\Z) \ar[r] & H^1(\Dbarre,\Q/\Z)^G \ar[r] & H^2(k,\Q/\Z) \ar[r] & H^2(D,\Q/\Z)\rlap{\text{.}}
}
$$
Si $D(k)\neq \emptyset$, le choix d'un point rationnel induit pour chaque~$i$
une r\'etraction de l'application de restriction $H^i(k,\Q/\Z)\to H^i(D,\Q/\Z)$.
En particulier la fl\`eche de droite de la suite exacte ci-dessus est alors
injective et l'on obtient une d\'ecomposition en somme directe
$H^1(D,\Q/\Z)=H^1(k,\Q/\Z)\oplus H^1(\Dbarre,\Q/\Z)^G$.  Si~$k$ est un corps
de type fini sur son sous-corps premier, le groupe $H^1(\Dbarre,\Q/\Z)^G$ est
fini (cf.~\cite[Th.~1~(bis)]{katzlang}).  C'est ce groupe fini que la premi\`ere
partie du lemme~\ref{lemme:cremona} d\'etermine dans un cas particulier.  Sous
les hypoth\`eses de ce lemme, l'assertion~(a) \'equivaut \`a dire que le groupe
$H^1(\Dbarre,\Q/\Z)^G$ est d'ordre~$3$.

(ii) Dans la situation de la proposition~\ref{prop:brauer-affine},
les groupes $\Br(\Xbarre)$ et $H^3(\Xbarre,\Gm)$ sont nuls puisque
$\Xbarre=X \otimes_\Q \Qbarre$ est une surface rationnelle.  Le th\'eor\`eme de
puret\'e pour le groupe de Brauer et la suite exacte de localisation
fournissent donc un isomorphisme de modules galoisiens
$\Br(\Ubarre)\simeq H^1(\Dbarre,\Q/\Z)$ (cf.~\cite[\textsection6]{brauerIII}
ou la suite exacte~(\ref{se:loc}) ci-dessous).  Ce dernier groupe est
isomorphe \`a $(\Q/\Z)^2$.  Compte tenu de la remarque~\ref{rq:HSC}~(i) et du
lemme~\ref{lemme:cremona}, on a ici $\Br(\Ubarre)^G \simeq \Z/3\Z$.
En particulier $\Br(\Ubarre)^G$ n'est pas trivial.
La proposition~\ref{prop:brauer-affine} entra\^ine pourtant que l'application
naturelle $\Br(U) \to \Br(\Ubarre)^G$ est nulle, puisque $\Br(\Xbarre)=0$.
\end{rmqs}

\begin{prop}
\label{prop:brauer-affine2}
Soit $a \in \Q^*$.
Notons $X' \subset \P^3_\Q$ la surface projective et lisse d'\'equation homog\`ene
$x^3+y^3+2z^3=at^3$ et $U' \subset \A^3_\Q$ la surface affine d'\'equation
$x^3+y^3+2z^3=a$, compl\'ementaire dans~$X'$ de la section hyperplane $t=0$.
L'image de l'application de restriction $\Br(X') \to \Br(U')$ est d'indice~$2$
dans $\Br(U')$.  La classe de l'alg\`ebre de quaternions
$$B'=(a(x+y+2z), -3(x+y+2z)(x+y)) \in \Br(\Q(U'))$$
appartient au sous-groupe $\Br(U') \subset \Br(\Q(U'))$ et n'appartient pas \`a
$\Br(X')$, de sorte qu'elle induit une d\'ecomposition en somme directe
$\Br(U')=\Br(X') \oplus \Z/2\Z$.
\end{prop}

Sous les hypoth\`eses de la proposition~\ref{prop:brauer-affine2}, si aucun
de~$a$, $2a$, $4a$ n'est un cube dans~$\Q$, le quotient $\Br(U')/\Br(\Q)$ est
donc d'ordre~$6$, engendr\'e par les images des classes~$A'$ et~$B'$ d\'efinies
respectivement dans les propositions~\ref{prop:proj-br}
et~\ref{prop:brauer-affine2}.

\begin{proof}
Le lemme~\ref{lemme:cremona2} ci-dessous jouera le m\^eme r\^ole ici que le
lemme~\ref{lemme:cremona} dans la preuve de la
proposition~\ref{prop:brauer-affine}:

\begin{lem}
\label{lemme:cremona2}
Notons $D' \subset \P^2_\Q$ la courbe plane d'\'equation $x^3+y^3+2z^3=0$.
Posons $P_0=[1:-1:0] \in D'(\Q)$ et $P_1=[1:0:-1/\sqrt[3]2\mkern2mu] \in D'(K)$,
o\`u $K=\Q(\sqrt[3]2\mkern2mu)$.

(a) Le noyau de la fl\`eche $H^1(D,\Q/\Z)\to H^1(\Q,\Q/\Z)$ d'\'evaluation
en~$P_0$ est isomorphe \`a $\Z/6\Z$.

(b) Le noyau de la fl\`eche $H^1(D,\Q/\Z)\to H^1(\Q,\Q/\Z) \times H^1(K,\Q/\Z)$
produit des fl\`eches d'\'evaluation aux points~$P_0$ et~$P_1$ est isomorphe \`a
$\Z/2\Z$.
\end{lem}

\begin{proof}
Munissons~$D'$ de la structure de courbe elliptique sur~$\Q$ ayant pour
\'el\'ement neutre le point~$P_0$.  Une \'equation de Weierstrass pour~$D'$ est
donn\'ee par
$\eta^2=\xi^3-1/27$ (\emph{via} $x=-1-3\eta$, $y=-1+3\eta$, $z=3\xi$).
Il~s'agit de la courbe 36A3 dans les tables de
Cremona~\cite{cremona}. Celles-ci nous apprennent que le sous-groupe de
torsion du groupe de Mordell--Weil de chaque courbe elliptique sur~$\Q$
isog\`ene \`a~$D'$ est cyclique, d'ordre~$2$ ou~$6$.  Par dualit\'e, le noyau de la
fl\`eche $H^1(D',\Q/\Z) \to H^1(\Q,\Q/\Z)$ d'\'evaluation en~$P_0$ s'identifie
donc \`a $\Hom(\mmu_6,D')$.

On v\'erifie sans peine que le groupe alg\'ebrique $\tors3D'$ s'ins\`ere dans une
suite exacte
$$
\xymatrix{
0 \ar[r] & \mmu_3 \ar[r] & \tors3D' \ar[r] & \Z/3\Z \ar[r] & 0\rlap{\text{,}}
}
$$
o\`u $\mmu_3 \subset \tors3D'$ d\'esigne l'intersection de~$D'$ avec l'hyperplan
d'\'equation $z=0$ dans~$\P^3_\Q$.
Il s'ensuit que $\Hom(\mmu_3,D')=\Hom(\mmu_3,\mmu_3)=\Z/3\Z$.
Par ailleurs, le groupe $\tors2D'(\Q)$ est d'ordre~$2$, de sorte que
$\Hom(\mmu_2,D')=\Z/2\Z$ et donc $\Hom(\mmu_6,D')=\Z/6\Z$.
Ainsi l'\'evaluation en~$P_0$ induit-elle une d\'ecomposition
$H^1(D',\Q/\Z)=H^1(\Q,\Q/\Z) \times \Z/6\Z$, le facteur $\Z/6\Z$ \'etant
engendr\'e par les classes $m_2,m_3 \in H^1(D',\Q/\Z)$ des suites exactes
\begin{align}
\label{se:m2}
\xymatrix{
0 \ar[r] & \Z/2\Z \ar[r] & D'/\mmu_2 \ar[r] & D' \ar[r] & 0
}
\end{align}
et
\begin{align}
\label{se:m3}
\xymatrix{
0 \ar[r] & \Z/3\Z \ar[r] & D'/\mmu_3 \ar[r] & D' \ar[r] & 0
}
\end{align}
induites respectivement par la multiplication par~$2$ et par~$3$ sur~$D'$,
o\`u $\mmu_2 \subset D'$ d\'esigne le sous-groupe $\tors2D'(\Q)$ de~$\tors2D'$.
Ceci \'etablit l'assertion~(a) du lemme.

\`A l'aide des formules de V\'elu~\cite{velu}, qui fournissent des \'equations pour
la courbe elliptique $D'/\mmu_3$ et pour l'isog\'enie $D'/\mmu_3 \to D'$
apparaissant dans~(\ref{se:m3}), on v\'erifie que la fibre de $D'/\mmu_3 \to D'$
en~$P_1$ est irr\'eductible.
L'image de~$P_1$ par l'application $D'(K) \to H^1(K,\Z/3\Z)$ bord
de~(\ref{se:m3}) est donc non nulle.  En revanche l'image de~$P_1$ par
l'application $D'(K)\to H^1(K,\Z/2\Z)$ bord de~(\ref{se:m2}) est nulle puisque
le point~$P_1$ est d'ordre~$3$ sur~$D'$.  En d'autres termes, nous constatons
que $m_3(P_1)\neq 0$ mais $m_2(P_1)=0$. L'assertion~(b) du lemme s'ensuit.
\end{proof}

\'Etablissons maintenant la proposition~\ref{prop:brauer-affine2}.  Comme dans
la preuve de la proposition~\ref{prop:brauer-affine}, la section hyperplane
$D' \subset X'$ d'\'equation $t=0$ est lisse et l'on dispose donc, par puret\'e,
d'une application r\'esidu $\partial_{D'} : \Br(U') \to H^1(D',\Q/\Z)$ dont le
noyau est~$\Br(X')$.  Soient $A \in \Br(U')$ et $m=\partial_{D'}(A)$.  Notons
$k=\Q(\sqrt[3]{4a}\mkern2mu)$, $K=\Q(\sqrt[3]2\mkern2mu)$ et
$K'=K(\sqrt[3]a\mkern2mu)$.  La droite $L'_0 \subset \P^3_k$ d'\'equations
$x+y=0$, $z=(\sqrt[3]a/\sqrt[3]2\mkern2mu)t$
est contenue dans $X' \otimes_\Q k$ et rencontre $D' \otimes_\Q k$ au seul
point $P_0 \otimes_\Q k$, avec intersection transverse en ce point.
Le m\^eme raisonnement que celui employ\'e dans la preuve de la
proposition~\ref{prop:brauer-affine} permet d'en d\'eduire que $m(P_0)=0$,
compte tenu de l'injectivit\'e de la fl\`eche de restriction
$H^1(\Q,\Q/\Z) \to H^1(k,\Q/\Z)$.  La droite $L_1 \subset \P^3_{K'}$
d'\'equations $x+\sqrt[3]2\mkern1mu{}z=0$, $y=\sqrt[3]a\mkern1mu{}t$
est contenue dans $X' \otimes_\Q {K'}$ et rencontre $D' \otimes_\Q K'$ au seul
point $P_1 \otimes_K K'$, avec intersection transverse en ce point.
L'extension $K'/K$ \'etant soit triviale, soit cubique et non cyclique, la
fl\`eche de restriction $H^1(K,\Q/\Z)\to H^1(K',\Q/\Z)$ est injective.  Comme
pr\'ec\'edemment il s'ensuit que $m(P_1)=0$.  Le lemme~\ref{lemme:cremona2}~(b)
montre maintenant que si $m\neq 0$, alors~$m$ est l'unique \'el\'ement de
$H^1(D,\Q/\Z)$ v\'erifiant $m(P_0)=0$ et $m(P_1)=0$.  Le quotient
$\Br(U')/\Br(X')$ est donc d'ordre au plus~$2$.

Pour conclure la preuve de la proposition~\ref{prop:brauer-affine2} il reste
seulement \`a v\'erifier que $B' \in \Br(U')$ et $B' \notin \Br(X')$.  De fa\c{c}on
\'evidente, la classe~$B'$ est non ramifi\'ee sur~$U'$ en dehors des courbes
irr\'eductibles $x+y=0$ et $x+y+2z=0$.  Que le r\'esidu de~$B'$ au point g\'en\'erique
de chacune de ces deux courbes soit nul r\'esulte de la formule du symbole
mod\'er\'e (cf.~\cite[Ex.~7.1.5, Prop.~7.5.1]{gilleszamuely}) et de l'identit\'e
polynomiale
\begin{align}
\label{eq:identite-polynomiale}
\begin{split}
\congru{4(x^3+y^3+2z^3)}{3(x+y)(x-y)^2}{x+y+2z}\rlap{\text{.}}
\end{split}
\end{align}
Quant au r\'esidu de~$B'$ au point g\'en\'erique de~$D'$, il est \'egal \`a la classe de
la fonction rationnelle $-3(x+y+2z)/(x+y)$ dans
$\Q(D')^*/\Q(D')^{*2}=H^1(\Q(D'),\Z/2\Z)$.  Le diviseur de cette fonction est
$2P-2P_0$ o\`u $P \in D'(\Q)$ a pour coordonn\'ees homog\`enes $[1:1:-1]$.  Comme
$P-P_0$ n'est pas un diviseur principal sur~$D'$, cette fonction n'est pas un
carr\'e et la proposition~\ref{prop:brauer-affine2} est \'etablie.
\end{proof}

\begin{rmq}
La fin de la d\'emonstration de la proposition~\ref{prop:brauer-affine2} montre
que le r\'esidu de~$B'$ au point g\'en\'erique de~$D'$ reste non nul m\^eme apr\`es
extension des scalaires de~$\Q$ \`a un corps alg\'ebriquement clos.  Par
cons\'equent~$B'$ est une classe \og{}transcendante\fg{} du groupe de Brauer de~$U'$,
contrairement \`a~$A'$.
\end{rmq}

\section{Sur les entiers}
\label{par:sur-les-entiers}

Nous sommes \`a pr\'esent en position d'\'etablir le

\begin{thm}
\label{thm:sur-les-entiers}
Soit $a \in \Z$ un entier non nul.

(a) Si~$a$ n'est pas de la forme $9n\pm 4$ pour un $n \in \Z$, il n'y
a pas d'obstruction de Brauer--Manin \`a l'existence d'une solution de
l'\'equation $x^3+y^3+z^3=a$ en entiers $x, y, z \in \Z$.

(b) Il n'y a pas d'obstruction de Brauer--Manin \`a l'existence d'une
solution de l'\'equation $x^3+y^3+2z^3=a$ en entiers $x, y, z \in \Z$.
\end{thm}

Rappelons que lorsque $a=9n\pm 4$ pour un $n \in \Z$, l'\'equation
$x^3+y^3+z^3=a$ n'admet pas de solution enti\`ere (elle n'admet m\^eme pas de
solution modulo~$9$).

Les r\'esultats des paragraphes~\ref{par:corps-arbitraire}
et~\ref{par:rationnels} fournissent des g\'en\'erateurs explicites du groupe de
Brauer des surfaces affines intervenant dans le th\'eor\`eme.  Afin de d\'eterminer
l'obstruction de Brauer--Manin enti\`ere associ\'ee \`a ces g\'en\'erateurs, nous nous
servirons d'un lemme g\'en\'eral permettant de relier le calcul des invariants
locaux \`a la g\'eom\'etrie modulo~$p$ des surfaces consid\'er\'ees, aux places de
mauvaise r\'eduction.  Nous \'enon\c{c}ons et d\'emontrons ce lemme
(lemme~\ref{lemme:ev-cone}) ainsi que deux compl\'ements
(lemmes~\ref{lemme:ev-cone-nonramifie} et~\ref{lemme:ev-cone-existence}) avant
d'entamer la preuve du th\'eor\`eme~\ref{thm:sur-les-entiers}.
Ces trois lemmes seront utilis\'es dans les preuves des
propositions~\ref{prop:evaluationA} et~\ref{prop:evalApBp} pour calculer
l'image de l'application qui \`a un $\Z_p$\nobreakdash-point~$M$ de la
surface affine d'\'equation $x^3+y^3+z^3=a$ (resp.~$x^3+y^3+2z^3=a$)
associe l'invariant $p$\nobreakdash-adique de~$A(M)$ (resp.~de~$A'(M)$ ou
de~$B'(M)$) lorsque~$p$ divise~$a$, o\`u~$A$, $A'$ et~$B'$ sont les classes
donn\'ees par les propositions~\ref{prop:proj-br}
et~\ref{prop:brauer-affine2}.

\begin{lem}
\label{lemme:ev-cone}
Soit~$K$ un corps $p$\nobreakdash-adique, d'anneau des entiers~$R$ et de corps
r\'esiduel~$k$.  Soient $n \geq 1$ et $V \subset \P^n_K$ une vari\'et\'e projective
et lisse.  Fixons un hyperplan $H \subset \P^n_K$ et posons $D = V \cap H$ et
$U = V \setminus D$.
Notons~$\sV$ et~$\sH$ les adh\'erences de~$V$ et de~$H$ dans $\P^n_R$.
Supposons le sch\'ema $\sD = \sV \cap \sH$ lisse sur~$R$ et supposons qu'il
existe $s \in \sV(k)$ n'appartenant pas \`a~$\sH$, tel que $\sV \otimes_R k$
soit le c\^one dans~$\P^n_k$ de sommet~$s$ et de base $\sD\otimes_R k$.  Notons
$\pi: \sV^0 \otimes_R k \to \sD \otimes_R k$ le morphisme de projection
depuis~$s$, o\`u $\sV^0=\sV \setminus \{s\}$.  Soient enfin $N \geq 1$ un entier
inversible dans~$R$ et $B \in \Br(U)$ une classe annul\'ee par~$N$, telle que la
condition~($*$) soit satisfaite:

\medskip
($*$) {\setbox0=\hbox{($*$)}\dimen0\textwidth\advance\dimen0-\parindent\advance\dimen0-\parindent\advance\dimen0-\wd0\parbox[t]{\dimen0}
{Le r\'esidu de~$B$ au point g\'en\'erique de~$D$ appartient \`a l'image de la fl\`eche de restriction $H^1(\sD,\Z/N\Z) \to H^1(D,\Z/N\Z)$.}}

\medskip
\noindent{}Il existe alors un unique $\gamma \in H^1(\sD \otimes_Rk,\Z/N\Z)$
tel que pour tout point g\'en\'erique~$\eta$ d'une composante irr\'eductible de~$\sV
\otimes_Rk$, le r\'esidu de~$B$ en~$\eta$ soit \'egal \`a
$(\pi^*\gamma)(\eta) \in H^1(\eta,\Z/N\Z)$. En outre, si l'on pose
$\sU^0=\sV^0 \setminus \sD$, alors pour tout $M \in U(K)$ qui se sp\'ecialise en
un point $m \in \sV(k)$ appartenant \`a~$\sU^0(k)$, l'image de $\gamma(\pi(m))$
par l'isomorphisme canonique $H^1(k,\Z/N\Z) \isoto \Z/N\Z$ co\"incide avec
l'invariant associ\'e \`a~$B(M)$ par la th\'eorie du corps de classes local.
\end{lem}

\begin{proof}
Notons $\sZ = \sV^0 \setminus U$ et $\sZ^0 = \sZ \setminus (\sD \otimes_Rk)$,
de sorte que $\sZ^0$ est la r\'eunion disjointe de~$D$ et de $\sU^0 \otimes_Rk$.
Les sch\'emas $\sZ^0$, $\sV^0$ et $\sD \otimes_Rk$ \'etant tous trois r\'eguliers,
la suite exacte longue du triple \cite[Ch.~III, Rem.~1.26]{milne}
$$
\xymatrix@C=2em{
\cdots \ar[r] & H^3_{\sZ}(\sV^0,\mmu_N) \ar[r] & H^3_{\sZ^0}(\sV^0 \setminus (\sD \otimes_Rk),\mmu_N) 
 \ar[r] & H^{4}_{\sD\otimes_Rk}(\sV^0,\mmu_N) \ar[r] & \cdots
}
$$
se r\'e\'ecrit, par puret\'e~\cite{fujiwara}, en
\begin{align}
\label{se:triple}
\xymatrix@R=0ex{
H^3_{\sZ}(\sV^0,\Lambda(1)) \ar[r] & H^1(D,\Lambda) \oplus H^1(\sU^0\otimes_Rk,\Lambda) \ar[r]^(.565){f\oplus g} & H^0(\sD \otimes_Rk, \Lambda(-1)) \rlap{\text{,}}
}
\end{align}
o\`u~$\Lambda(r)$ d\'esigne le faisceau $\Lambda=\Z/N\Z$ tordu~$r$ fois \`a la Tate.
De plus, toujours par puret\'e, les fl\`eches~$f$ et~$g$ s'inscrivent dans des
suites exactes
\begin{align}
\label{se:purete-D}
\xymatrix{
0 \ar[r] & H^1(\sD,\Lambda) \ar[r] & H^1(D,\Lambda) \ar[r]^(.37)f & H^0(\sD\otimes_Rk,\Lambda(-1))
}
\end{align}
et
\begin{align}
\label{se:purete-sU}
\xymatrix{
0 \ar[r] & H^1(\sV^0\otimes_Rk,\Lambda)\ar[r]&H^1(\sU^0\otimes_Rk,\Lambda)\ar[r]^(.45)g&H^0(\sD\otimes_Rk,\Lambda(-1))\rlap{\text{.}}
}
\end{align}
Notons
$\alpha \oplus \beta \in H^1(D,\Lambda) \oplus H^1(\sU^0\otimes_Rk,\Lambda)$
l'image, par la compos\'ee de l'application naturelle
$H^2(U,\Lambda(1)) \to H^3_{\sZ}(\sV^0,\Lambda(1))$
et de la fl\`eche de gauche de~(\ref{se:triple}), d'un rel\`evement de~$B$ dans
$H^2(U,\Lambda(1))$.  Cette image ne d\'epend pas du rel\`evement choisi; plus
pr\'ecis\'ement, les restrictions de~$\alpha$ et de~$\beta$ aux points g\'en\'eriques
de~$D$ et des composantes irr\'eductibles de $\sV\otimes_Rk$ ne sont autres que
les r\'esidus de~$B$ en ces points.  Au vu de la suite~(\ref{se:purete-D}), la
condition~($*$) entra\^ine que $f(\alpha)=0$.  D'autre part,
comme~(\ref{se:triple}) est un complexe, on a $f(\alpha)+g(\beta)=0$.  Ainsi
$g(\beta)=0$, ce qui, d'apr\`es la suite exacte~(\ref{se:purete-sU}), signifie
que $\beta \in H^1(\sV^0\otimes_Rk,\Lambda)$.

La fl\`eche $\pi^*:H^1(\sD\otimes_Rk,\Lambda) \to H^1(\sV^0\otimes_Rk,\Lambda)$
est un isomorphisme puisque~$\pi$ est le morphisme structural d'un fibr\'e en
droites (cf.~\cite[Ch.~VI, Cor.~4.20]{milne}).  Il~existe donc un unique
$\gamma \in H^1(\sD\otimes_Rk,\Lambda)$ tel que $\beta=\pi^*\gamma$.

Soit maintenant $M \in U(K)$ se sp\'ecialisant en un point $m \in \sU^0(k)$.
Notons $\sM \subset \sU^0$ l'adh\'erence de~$M$ dans~$\sV$, de sorte que
$\sM\otimes_Rk=m$.  Comme~$\sU^0$ et~$\sM$ sont lisses sur~$R$, on dispose par
puret\'e d'un carr\'e commutatif
$$
\xymatrix{
H^2(U,\mmu_N) \ar[r] \ar[d] & H^1(\sU^0 \otimes_Rk, \Z/N\Z) \ar[d]^e\\
H^2(M,\mmu_N) \ar[r] & H^1(\sM \otimes_Rk, \Z/N\Z)\rlap{\text{,}}
}
$$
dont les fl\`eches verticales sont les applications de restriction et dont les
fl\`eches horizontales sont les applications r\'esidu.  La commutativit\'e de ce
carr\'e implique que l'image de $e(\beta)$ par l'isomorphisme canonique
$H^1(\sM\otimes_Rk,\Z/N\Z)=H^1(k,\Z/N\Z)\isoto \Z/N\Z$ est \'egale \`a l'invariant
de~$B(M)$ (cf.~\cite[p.~94]{brauerIII}).  Comme par ailleurs
$e(\beta)=\beta(m)=(\pi^*\gamma)(m)=\gamma(\pi(m))$, ceci termine la
d\'emonstration du lemme.
\end{proof}

\begin{lem}
\label{lemme:ev-cone-nonramifie}
Sous les hypoth\`eses du lemme~\ref{lemme:ev-cone}, si la classe~$B$ appartient
au sous-groupe $\Br(V) \subset \Br(U)$ (de sorte que la condition~($*$) est
trivialement satisfaite) et si $B|_D\in\Br(D)$ d\'esigne la restriction de~$B$
\`a~$D$, alors~$\gamma$ est \'egal au r\'esidu de~$B|_D$ au point g\'en\'erique
de $\sD \otimes_Rk$.
\end{lem}

\begin{proof}
Les sch\'emas~$\sV^0$ et~$\sD$ \'etant lisses sur~$R$, le th\'eor\`eme de puret\'e
absolu~\cite{fujiwara} fournit les fl\`eches horizontales du carr\'e commutatif
$$
\xymatrix{
H^2(V,\mmu_N) \ar[r] \ar[d] & H^1(\sV^0 \otimes_Rk, \Z/N\Z) \ar[d]_{i^*} \\
H^2(D,\mmu_N) \ar[r] & \ar@/_1.1pc/@{.>}[u]_{\pi^*} H^1(\sD\otimes_Rk,\Z/N\Z)\rlap{\text{,}}
}
$$
o\`u~$i$ d\'esigne l'immersion ferm\'ee $\sD\otimes_Rk\to\sV^0\otimes_Rk$.
De la commutativit\'e de ce carr\'e s'ensuit
que le r\'esidu de $B|_D$ au point g\'en\'erique de $\sD\otimes_Rk$ est \'egal
\`a~$i^*\beta$. Or
$i^*\beta=i^*\pi^*\gamma=(\pi \circ i)^* \gamma=\gamma$, d'o\`u le lemme.
\end{proof}

\begin{lem}
\label{lemme:ev-cone-existence}
Pla\c{c}ons-nous dans la situation du lemme~\ref{lemme:ev-cone} et fixons une
famille finie $\uplet{B_1}{B_\ell} \in \Br(U)$ de classes annul\'ees par~$N$,
v\'erifiant chacune la condition~($*$).  Pour tout~$i$, notons~$\gamma_i$ la
classe associ\'ee \`a~$B_i$ par le lemme~\ref{lemme:ev-cone} et~$N_i$ l'ordre
de~$\gamma_i$.  Supposons que~$D$ soit une courbe elliptique sur~$K$, que les
entiers~$N_i$ soient premiers entre eux deux \`a deux et que pour tout~$i$,
l'image de~$\gamma_i$ dans $H^1(\sD\otimes_R\kbar,\Z/N\Z)$ soit d'ordre~$N_i$,
o\`u~$\kbar$ d\'esigne une cl\^oture alg\'ebrique de~$k$.  Alors l'application
$$\sU^0(R) \longrightarrow (\Z/N\Z)^\ell$$
qui envoie $M \in \sU^0(R) \subset U(K)$ sur la famille des invariants associ\'es
aux $B_i(M)$ par la th\'eorie du corps de classes local a pour image
$\prod_{i=1}^\ell \Z/N_i\Z$.
\end{lem}

\begin{proof}
Pour tout~$i$, la classe $\gamma_i \in H^1(\sD\otimes_Rk,\Z/N\Z)$ est
repr\'esent\'ee par un rev\^etement cyclique $E_i \to \sD\otimes_Rk$ de degr\'e~$N_i$.
Ce rev\^etement est g\'eom\'etriquement connexe sur~$k$ puisque l'image
de~$\gamma_i$ dans $H^1(\sD\otimes_R\kbar,\Z/N\Z)$ est encore d'ordre~$N_i$.
De plus,
comme le corps~$k$ est fini,
il poss\`ede un point rationnel,
d'apr\`es un th\'eor\`eme de F.~K.~Schmidt
(cf.~\cite[Ch.~III, \textsection2.3]{cohomologiegaloisienne}).
Il s'agit donc d'une courbe elliptique s'inscrivant dans une suite exacte
\begin{align}
\label{se:gamma}
\xymatrix{
0 \ar[r] & \Z/N_i\Z \ar[r] & E_i \ar[r] & \sD \otimes_Rk \ar[r] & 0\rlap{\text{.}}
}
\end{align}
De ce point de vue, l'application $\sD(k) \to H^1(k,\Z/N_i\Z)$,
$m_0 \mapsto \gamma_i(m_0)$ n'est autre que le bord de~(\ref{se:gamma}); c'est
donc un morphisme de groupes.  Celui-ci est surjectif puisque l'on a
$H^1(k,E_i)=0$ d'apr\`es le th\'eor\`eme de F.~K.~Schmidt.  Comme
les~$N_i$ sont premiers entre eux deux \`a deux, il s'ensuit que l'application
$\sD(k) \to \prod_{i=1}^\ell H^1(k,\Z/N_i\Z)$,
$m_0 \mapsto (\gamma_i(m_0))_{1 \leq i \leq \ell}$ est elle aussi surjective.
D'autre part la fl\`eche de sp\'ecialisation $\sU^0(R) \to \sU^0(k)$ est
surjective (lemme de Hensel) et le morphisme
$\sU^0 \otimes_Rk \to \sD\otimes_Rk$ induit par~$\pi$ est surjectif sur les
$k$\nobreakdash-points (ses fibres \'etant isomorphes \`a~$\Gm$).
La compos\'ee de ces trois surjections, de l'isomorphisme canonique
$\prod_{i=1}^\ell H^1(k,\Z/N_i\Z) \isoto \prod_{i=1}^\ell \Z/N_i\Z$
et de l'inclusion $\prod_{i=1}^\ell \Z/N_i\Z \subset (\Z/N\Z)^\ell$ est une
application $\sU^0(R) \longrightarrow (\Z/N\Z)^\ell$ dont le
lemme~\ref{lemme:ev-cone} montre qu'elle co\"incide avec celle apparaissant dans
l'\'enonc\'e du lemme~\ref{lemme:ev-cone-existence}, ce qui conclut la
d\'emonstration.
\end{proof}

\begin{proof}[D\'emonstration du th\'eor\`eme~\ref{thm:sur-les-entiers}]
Si $a \in \Z$, posons $\sU_a=\Spec(\Z[x,y,z]/(x^3+y^3+z^3-a))$
et $\sU'_a=\Spec(\Z[x,y,z]/(x^3+y^3+2z^3-a))$.

\begin{lem}
\label{lemme:points-locaux}
Les ensembles $\sU_a(\R)$, $\sU'_a(\R)$, $\sU_a(\Z_p)$ et $\sU'_a(\Z_p)$
sont non vides pour tout~$a$ et tout~$p$,
\`a l'exception de $\sU_a(\Z_3)$ dans le cas o\`u $\congru{a}{\pm 4}{9}$.
\end{lem}

\begin{proof}
Pour tout $p\neq 3$, le sch\'ema $\sU_a$ admet un $\Fp$\nobreakdash-point
lisse. En effet, si~$p$ divise~$a$, le point $(1,-1,0)$ convient; dans le cas
contraire, tout $\Fp$\nobreakdash-point est lisse et il est bien connu que
$\sU_a(\Fp)\neq\emptyset$ (cf.~\cite[Th.~1]{tornheim}).  Il s'ensuit que
$\sU_a(\Z_p)\neq \emptyset$ pour tout $p \neq 3$.  Pour $p \notin\{2,3\}$,
l'ensemble $\sU'_a(\Z_p)$ est non vide en vertu de
l'identit\'e~(\ref{eq:identite-6a}) appliqu\'ee \`a~$a/6$.  Que~$\sU_a(\R)$ et
$\sU'_a(\R)$ soient non vides est \'evident.  Enfin, la non vacuit\'e de
$\sU_a(\Z_3)$ lorsque $\pascongru{a}{\pm 4}{9}$ et celle de~$\sU'_a(\Z_2)$ et
de $\sU'_a(\Z_3)$ en toute g\'en\'eralit\'e se ram\`enent \emph{via} le lemme de
Hensel \`a un calcul fini que nous ne d\'etaillons pas ici.
\end{proof}

Quels que soient les entiers~$a$ et~$n$, il r\'esulte de l'existence d'un
morphisme $\sU_a \to \sU_{an^3}$ (resp.~$\sU'_a \to \sU'_{an^3}$) que s'il n'y
a pas d'obstruction de Brauer--Manin \`a l'existence d'une solution de
l'\'equation $x^3+y^3+z^3=a$ (resp.~$x^3+y^3+2z^3=a$) en entiers
$x, y, z \in \Z$, il n'y a pas davantage d'obstruction de Brauer--Manin
\`a l'existence d'une solution de l'\'equation $x^3+y^3+z^3=an^3$
(resp.~$x^3+y^3+2z^3=an^3$) en entiers $x, y, z \in \Z$.  En vue d'\'etablir le
th\'eor\`eme, et compte tenu que $\congru{a}{\pm 4}{9}$ entra\^ine
$\congru{an^3}{\pm 4}{9}$ si~$n$ est premier \`a~$3$, il est donc loisible de
supposer que pour tout nombre premier $p \neq 3$, la valuation
$p$\nobreakdash-adique de~$a$ v\'erifie $v_p(a) \in \{0,1,2\}$.  Nous fixons
dor\'enavant un entier $a>0$ satisfaisant cette hypoth\`ese, et posons $\sU=\sU_a$
et $\sU'=\sU'_a$.  Nous reprenons les notations $U$, $U'$, $X$, $X'$, $D$,
$D'$, $A$, $A'$, $B'$ introduites dans les paragraphes pr\'ec\'edents et notons de
plus $\sX$, $\sX'$ les sous-sch\'emas ferm\'es de~$\P^3_\Z$ d'\'equations homog\`enes
respectives $x^3+y^3+z^3=at^3$ et $x^3+y^3+2z^3=at^3$, et $\sD$, $\sD'$ les
sous-sch\'emas ferm\'es de~$\sX$ et de~$\sX'$ d\'efinis par l'\'equation $t=0$.  Ainsi
$U = \sU \otimes_\Z \Q$, $U'=\sU' \otimes_\Z \Q$, $\sU=\sX \setminus \sD$ et
$\sU'=\sX' \setminus \sD'$.

Consid\'erons d'abord l'obstruction de Brauer--Manin enti\`ere sur~$\sU$.  Si~$a$
est une puissance de~$3$ alors $\sU(\Z)\neq\emptyset$ de mani\`ere \'evidente et
il n'y a rien \`a d\'emontrer.  Sinon, il existe un nombre premier $q \neq 3$ tel
que $v_q(a) \in \{1,2\}$.  Dans ce cas~$a$ n'est pas un cube et les
propositions~\ref{prop:proj-br} et~\ref{prop:brauer-affine} entra\^inent que le
quotient $\Br(U)/\Br(\Q)$ est d'ordre~$3$, engendr\'e par la classe de~$A$.
Soit alors $(M_p)_{p \in \Omega} \in \prod_{p \in \Omega} \sU(\Z_p)$ une
famille de points locaux, o\`u~$\Omega$ d\'esigne l'ensemble des places de~$\Q$
(on convient que $\Z_\infty=\R$); l'existence de $(M_p)_{p \in \Omega}$ est
assur\'ee par le lemme~\ref{lemme:points-locaux}.  D'apr\`es la
proposition~\ref{prop:evaluationA} ci-dessous, il existe $M'_q \in \sU(\Z_q)$
tel que $\inv_q A(M'_q) = -\sum_{p \in \Omega\setminus\{q\}}\inv_p A(M_p)$.
Ainsi, en posant $M'_p=M_p$ pour $p \neq q$ on obtient une famille
$(M'_p)_{p\in \Omega}\in\prod_{p\in\Omega}\sU(\Z_p)$ orthogonale \`a~$\Br(U)$
pour l'accouplement de Brauer--Manin, ce qui \'etablit l'assertion~(a) du
th\'eor\`eme.

\begin{prop}
\label{prop:evaluationA}
Pour $p \neq 3$ tel que $v_p(a) \in \{1,2\}$,
l'application $\sU(\Z_p) \to \Z/3\Z$ qui \`a $M \in \sU(\Z_p) \subset U(\Q_p)$
associe l'invariant $p$\nobreakdash-adique de $A(M) \in \Br(\Q_p)$ est
surjective.
\end{prop}

\begin{proof}
Comme $p \neq 3$, le sch\'ema $\sD \otimes_\Z\Z_p$ est lisse sur~$\Z_p$.
Par puret\'e on dispose donc d'une fl\`eche r\'esidu
$\tors{3}\Br(D\otimes_\Q\Q_p)\to H^1(\sD\otimes_\Z\Fp,\Z/3\Z)$.
Soit $\gamma$ l'image, par cette fl\`eche, de la restriction de $A\in\Br(X)$ \`a
$D \otimes_\Q\Q_p$.

\begin{lem}
\label{lemme:residu-geom-connexe}
L'image de~$\gamma$ dans $H^1(\sD\otimes_\Z\Fbarp,\Z/3\Z)$ est non nulle
(o\`u $\Fbarp$ d\'esigne une cl\^oture alg\'ebrique de~$\Fp$).
\end{lem}

\begin{proof}
Apr\`es extension des scalaires de~$\Q$ \`a~$\Q(j)$, la classe~$A$ s'\'ecrit
$(1+\sigma)\left(((x+jy)/(x+y),a)_j\right)=((x+jy)/(x+y),a)_j
- ((x+j^2y)/(x+y),a)_j
= ((x+jy)/(x+j^2y),a)_j$
d'apr\`es l'\'egalit\'e~(\ref{eq:action}).
Fixons un plongement de~$\Q(j)$ dans~$\Q_p^\nr$ et notons~$L$ le corps des
fonctions de $\sD\otimes_\Z\Fbarp$.  La fonction rationnelle
$f=(x+jy)/(x+j^2y)$ sur $\sD\otimes_\Z\Z_p^\nr$ est inversible au point
g\'en\'erique de $\sD\otimes_\Z\Fbarp$.  Par cons\'equent le r\'esidu, en ce point, de
la restriction de~$A$ \`a $D\otimes_\Q\Q_p^\nr$ est \'egal \`a la classe de
$f^{v_p(a)}$ dans $L^*/L^{*3}= H^1(L,\mmu_3)= H^1(L,\Z/3\Z)$ (o\`u $\Z/3\Z$ est
identifi\'e \`a $\mmu_3$ au-dessus de~$\Fbarp$ \`a l'aide du choix pr\'ec\'edemment fix\'e
de~$j$ dans $\Q_p^\nr$).  Il suffit donc, pour conclure, de v\'erifier que~$f$
n'est pas un cube dans~$L$.  Or le diviseur de~$f$ sur $\sD\otimes_\Z\Fbarp$
s'\'ecrit $3P-3Q$ avec $P,Q \in \sD(\Fbarp)$ distincts, et deux points distincts
sur une courbe elliptique ne sont jamais lin\'eairement \'equivalents.
\end{proof}

Compte tenu du lemme~\ref{lemme:residu-geom-connexe}, il suffit \`a pr\'esent
d'appliquer les lemmes~\ref{lemme:ev-cone-nonramifie}
et~\ref{lemme:ev-cone-existence} (avec $\ell=1$) pour conclure la preuve de la
proposition~\ref{prop:evaluationA}.
\end{proof}

Passons maintenant \`a la d\'emonstration de l'assertion~(b) du
th\'eor\`eme~\ref{thm:sur-les-entiers}.  Si~$a$ est une puissance de~$2$ ou une
puissance de~$3$, alors $\sU(\Z)\neq\emptyset$.  Il en va de m\^eme si~$a$ est
divisible par~$6$, compte tenu de~(\ref{eq:identite-6a}).  On peut donc
supposer que~$a$ ne s'\'ecrit pas sous la forme $\pm 2^r3^s$ avec $r,s \in \N$.
Dans ce cas ni~$a$, ni~$2a$, ni~$4a$ ne sont des cubes.  D'apr\`es les
propositions~\ref{prop:proj-br} et~\ref{prop:brauer-affine2}, le quotient
$\Br(U')/\Br(\Q)$ est donc d'ordre~$6$ et est engendr\'e par les classes de~$A'$
et de~$B'$.

\begin{prop}
\label{prop:evalApBp}
Soit $p \neq 2,3$.  L'application $\sU'(\Z_p)\to \Z/3\Z \times \Z/2\Z$ qui \`a
$M \in \sU'(\Z_p) \subset U(\Q_p)$ associe le couple form\'e des invariants
$p$\nobreakdash-adiques de $A'(M)$ et de $B'(M)$ est surjective si
$v_p(a)=1$. Si $v_p(a)=2$, son image est \'egale \`a $\Z/3\Z \times \{0\}$.
\end{prop}

\begin{proof}
Notons $\gamma_{A'} \in H^1(\sD'\otimes_\Z\Fp,\Z/3\Z)$ le r\'esidu, le long de
$\sD'\otimes_\Z\Fp$, de la restriction de $A' \in \Br(X')$ \`a
$D' \otimes_\Q\Q_p$.  Il r\'esulte du lemme~\ref{lemme:residu-geom-connexe},
moyennant un changement de variables \'evident, que l'image de~$\gamma_{A'}$
dans $H^1(\sD'\otimes_\Z\Fbarp,\Z/3\Z)$ est non nulle.

Le r\'esidu de~$B'$ au point g\'en\'erique de~$D'$ est la classe de la fonction
rationnelle $g=-3(x+y+2z)/(x+y)$ dans $\Q(D')^*/\Q(D')^{*2}$.  On voit gr\^ace
\`a~(\ref{eq:identite-polynomiale}) que le diviseur de~$g$ sur le sch\'ema
$\sD' \otimes_\Z \Z_p$ est un double.  Par cons\'equent, le rev\^etement double de
$\sD' \otimes_\Z \Z_p$ obtenu en extrayant une racine carr\'ee de~$g$ est \'etale
sur $\sD'\otimes_\Z\Z_p$; la classe~$B'$ v\'erifie donc la condition~($*$) du
lemme~\ref{lemme:ev-cone}. Notons~$\gamma_{B'}$ l'\'el\'ement de
$H^1(\sD'\otimes_\Z\Fp,\Z/2\Z)$ que celui-ci associe \`a~$B'$.

Supposons d'abord que $v_p(a)=1$.  Soient~$K$ le corps des fonctions de
$\sD'\otimes_\Z\Fbarp$ et~$L$ celui de $\sX'\otimes_\Z\Fbarp$, de sorte
que~$L$ est isomorphe \`a~$K(t)$.  Compte tenu que $v_p(a)=1$, le r\'esidu de~$B'$
au point g\'en\'erique~$\eta$ de $\sX' \otimes_\Z\Fbarp$ est la classe dans
$L^*/L^{*2}$ de la fonction rationnelle $-3(x+y+2z)/(x+y) \in K^*$.  Le
diviseur de cette fonction sur $\sD' \otimes_\Z\Fbarp$ est $2P-2P_0$ o\`u
$P,P_0 \in \sD'(\Fbarp)$ ont pour coordonn\'ees homog\`enes respectives $[1:1:-1]$
et $[1:-1:0]$.  Comme $P-P_0$ n'est pas un diviseur principal sur
$\sD'\otimes_\Z\Fp$, il s'ensuit que $-3(x+y+2z)/(x+y)$ n'est pas un carr\'e
dans~$K$, ni \emph{a fortiori} dans~$L$.  En conclusion, le r\'esidu de~$B'$
en~$\eta$ n'est pas nul et l'image de~$\gamma_{B'}$ dans
$H^1(\sD'\otimes_\Z\Fbarp,\Z/2\Z)$ est donc non nulle \'egalement.  Le
lemme~\ref{lemme:ev-cone-existence} appliqu\'e \`a~$\gamma_{A'}$ et~$\gamma_{B'}$
(avec $\ell=2$ et $N=6$) entra\^ine maintenant la surjectivit\'e de l'application
apparaissant dans l'\'enonc\'e de la proposition~\ref{prop:evalApBp}.

Supposons maintenant que $v_p(a)=2$.  Dans ce cas le r\'esidu de~$B'$ en~$\eta$
est nul.  Par cons\'equent $\gamma_{B'}=0$ et le lemme~\ref{lemme:ev-cone}
entra\^ine donc que l'application $\sU'^0(\Z_p)\to \Z/2\Z$ qui \`a~$M$ associe
l'invariant $p$\nobreakdash-adique de~$B'(M)$ est nulle, o\`u $\sU'^0$ d\'esigne
le compl\'ementaire, dans $\sU'$, du point de coordonn\'ees $(x,y,z)=(0,0,0)$ dans
$\sU'\otimes_\Z\Fp$.  Or on constate tout de suite que
$\sU'^0(\Z_p)=\sU'(\Z_p)$ puisque $v_p(a)<3$.  Il~suffit donc pour conclure de
v\'erifier que l'application $\sU'(\Z_p)\to\Z/3\Z$ d'\'evaluation de~$A'$ est
surjective; mais ceci r\'esulte du lemme~\ref{lemme:ev-cone-nonramifie}.
\end{proof}

\begin{prop}
\label{prop:evalBpdyadique}
L'application $\sU'(\Z_2) \to \Z/2\Z$ qui \`a~$M$ associe l'invariant
\mbox{dyadique} de $B'(M)$ est surjective.
\end{prop}

\begin{proof}
Toute unit\'e dyadique \'etant un cube, il~existe $v \in \N$ et $r \in \Z_2^*$
tels que $a=2^vr^3$. On a m\^eme $v \in \{0,1,2\}$ d'apr\`es les r\'eductions
effectu\'ees au d\'ebut de la d\'emonstration du th\'eor\`eme~\ref{thm:sur-les-entiers}.

Comme le morphisme canonique $\sU' \otimes_\Z\Z_2 \to \sU'_{2^v}
\otimes_\Z\Z_2$ consistant \`a diviser $x$, $y$ et~$z$ par~$r$ induit une
bijection $\sU'(\Z_2) \to \sU'_{2^v}(\Z_2)$ et comme l'image r\'eciproque de la
classe
$B'_{2^v}=(2^v(x+y+2z), -3(x+y+2z)(x+y)) \in \Br(\sU'_{2^v} \otimes_\Z \Q_2)$
par ce morphisme est \'egale \`a~$B'$, on peut supposer, quitte \`a remplacer~$a$
par~$2^v$, que $a \in \{1,2,4\}$.

Notons~$t$ la racine cubique de~$3$ dans $\Z_2^*$ et $\ev$ l'application
apparaissant dans l'\'enonc\'e du lemme.  Remarquant que $\congru{t}{27}{32}$,
il est ais\'e de v\'erifier les trois assertions suivantes, qui terminent la
d\'emonstration:
si $a=1$, alors $\ev(1,0,0)=0$ et $\ev(2,-1,-t)=1$;
si $a=2$, alors $\ev(1,1,1)=0$ et $\ev(t,1,-1)=1$;
si $a=4$, alors $\ev(0,0,1)=0$ et $\ev(t,-1,1)=1$.
\end{proof}

Nous pouvons maintenant conclure la preuve du
th\'eor\`eme~\ref{thm:sur-les-entiers}.  Rappelons que nous avons suppos\'e d'une
part que $v_p(a)\in \{0,1,2\}$ pour tout $p \neq 3$ et d'autre part que~$a$ ne
s'\'ecrit pas sous la forme $\pm 2^r3^s$.  Il existe donc un nombre premier
$q \neq 2,3$ tel que $v_q(a) \in \{1,2\}$.
Fixons alors, avec l'aide du lemme~\ref{lemme:points-locaux}, une famille de
points locaux $(M_p)_{p \in \Omega} \in \prod_{p \in \Omega} \sU'(\Z_p)$.  Si
$v_q(a)=1$, il existe, d'apr\`es la proposition~\ref{prop:evalApBp}, un point
$M'_q \in \sU'(\Z_q)$ tel que
$\inv_q A'(M'_q)=-\sum_{p \in \Omega\setminus\{q\}}\inv_p A'(M_p)$
et
$\inv_q B'(M'_q)=-\sum_{p \in \Omega\setminus\{q\}}\inv_p B'(M_p)$.
En posant $M'_p=M_p$ pour $p \neq q$ on obtient ainsi une famille
$(M'_p)_{p \in \Omega}\in\prod_{p \in \Omega}\sU'(\Z_p)$ orthogonale \`a la fois
\`a~$A'$ et \`a~$B'$, et donc \`a $\Br(U')$, pour l'accouplement de Brauer--Manin.
Si $v_q(a)=2$, il existe, d'apr\`es la proposition~\ref{prop:evalBpdyadique}, un
point $M'_2 \in \sU'(\Z_2)$ tel que
$\inv_2 B'(M'_2)=-\sum_{p \in \Omega\setminus\{2,q\}}\inv_p B'(M_p)$.
La proposition~\ref{prop:evalApBp} fournit ensuite un point
$M'_q \in \sU'(\Z_q)$ tel que
$\inv_q A'(M'_q)=-\inv_2 A'(M'_2) -\sum_{p \in \Omega\setminus\{2,q\}} \inv_p A'(M_p)$
et $\inv_q B'(M'_q)=0$.  Posant $M'_p=M_p$ pour $p \in \Omega \setminus \{2,q\}$,
la famille $(M'_p)_{p \in \Omega}$ est alors orthogonale \`a $\Br(U')$ pour
l'accouplement de Brauer--Manin.
\end{proof}

\begin{rmq}
\label{rq:evalB}
Nous avons vu avec la proposition~\ref{prop:evalApBp} que pour $p\geq 5$
premier, l'application $\sU'(\Z_p)\to \Z/2\Z$ qui \`a~$M$ associe l'invariant
$p$\nobreakdash-adique de~$B'(M)$ est nulle si~$v_p(a)$ est pair et est
surjective dans le cas contraire.  On peut v\'erifier que pour $p=3$ cette
application est nulle si $v_3(a) \leq 1$, surjective sinon,
et que pour $p=\infty$, elle est toujours nulle\footnote{Pour $p=3$ dans le
cas o\`u $v_3(a)=0$ et pour $p=\infty$, ces deux assertions r\'esultent
imm\'ediatement de l'identit\'e
$4(x^3+y^3+2z^3) = 3(x+y)(x-y)^2 + (x+y+2z)((x+y+2z)^2 - 6(x+y)z)$.}.
Ainsi, dans le cas o\`u~$a$ est un carr\'e premier \`a~$3$, la d\'emonstration du
th\'eor\`eme~\ref{thm:sur-les-entiers} passe-t-elle n\'ecessairement par des
consid\'erations dyadiques.
\end{rmq}

\section{Remarques, exemples et questions}
\label{par:remarques-questions}

\subsection{Groupes de Brauer}

Au paragraphe~\ref{par:rationnels} nous avons d\'etermin\'e, dans deux cas
particuliers, le conoyau de l'application de restriction $\Br(X) \to \Br(U)$
o\`u $X \subset \P^3_\Q$ est une surface cubique lisse et~$U$ est le
compl\'ementaire d'une section hyperplane lisse $D \subset X$
(propositions~\ref{prop:brauer-affine} et~\ref{prop:brauer-affine2}).  Les
d\'emonstrations s'appuyaient sur des informations de nature arithm\'etique (\'etude
des courbes elliptiques sur~$\Q$ isog\`enes \`a~$D$) et sur l'observation
alg\'ebrique suivante: si $L \subset \P^3_k$ est une droite contenue dans~$X$,
d\'efinie sur une extension $k/\Q$, alors le r\'esidu de toute classe de $\Br(U)$
le long de~$D$ s'annule au point d'intersection de~$L$ avec~$D$
puisque $L \cap U$ est isomorphe \`a~$\A^1_k$ et que $\Br(\A^1_k)=\Br(k)$.
Nous expliquons et g\'en\'eralisons cette observation dans le
lemme~\ref{lemme:controleresidu} ci-dessous; nous appliquons ensuite ce lemme dans
diverses situations.

Soient~$X$ une vari\'et\'e propre et lisse sur un corps~$k$ de caract\'eristique~$0$
et $D \subset X$ une sous-vari\'et\'e ferm\'ee, lisse, purement de codimension~$1$.
Notons $U=X \setminus D$.  D'apr\`es le th\'eor\`eme de puret\'e pour le groupe de
Brauer, la suite exacte de localisation en cohomologie \'etale s'\'ecrit
\begin{align}
\label{se:loc}
\xymatrix{
0 \ar[r] & \Br(X) \ar[r] & \Br(U) \ar[r] & H^1(D,\Q/\Z) \ar[r]^{\theta} & H^3(X,\Gm) \rlap{\text{.}}
}
\end{align}
(cf.~\cite[\textsection6]{brauerIII}).
Afin de d\'eterminer le groupe $\Br(U)$ il est n\'ecessaire de contr\^oler le noyau
de l'application~$\theta$ apparaissant dans~(\ref{se:loc}).  Pour ce faire,
nous consid\'ererons des courbes auxiliaires $C \subset X$ propres et lisses
telles que le sch\'ema $C \cap D$ soit \'etale sur~$k$.  Si~$C$ est une telle
courbe, notons $\sigma_C : H^1(D,\Q/\Z) \to H^1(k,\Q/\Z)$ l'application
d\'efinie par $\sigma_C(m)=\sum_{P \in C \cap D} \Cores_{k(P)/k}(m(P))$.

\begin{lem}
\label{lemme:controleresidu}
On a $\Ker(\theta) \subset \Ker(\sigma_C)$ pour toute courbe irr\'eductible
$C \subset X$ propre et lisse telle que le sch\'ema $C \cap D$ soit \'etale
sur~$k$.
\end{lem}

Dans le cas o\`u~$C$ et~$D$ se rencontrent en un unique point, on retrouve
l'observation signal\'ee pr\'ec\'edemment.

Pour prouver le lemme~\ref{lemme:controleresidu}, nous allons \'etablir l'\'enonc\'e
plus fort suivant:

\begin{lem}
\label{lemme:controleresidufort}
Soit $C \subset X$ une courbe irr\'eductible, propre et lisse sur~$k$, telle que le sch\'ema $C \cap D$ soit \'etale sur~$k$.
Notant~$\Cbarre$ la courbe obtenue par
extension des scalaires \`a une cl\^oture alg\'ebrique de~$k$,
on dispose d'un diagramme commutatif
$$
\xymatrix{
H^1(D,\Q/\Z) \ar[d]^{\sigma_C}\ar[r]^\theta & H^3(X,\Gm) \ar[r] & H^3(C,\Gm) \ar[r] & H^2(k,\Pic(\Cbarre)) \ar[d] \\
H^1(k,\Q/\Z) \ar[rrr]^\sim &&& H^2(k,\Z) \rlap{\text{,}}
}
$$
o\`u
la fl\`eche verticale de droite est induite par l'application degr\'e $\Pic(\Cbarre)\to\Z$ et
les fl\`eches horizontales du haut sont
la restriction $H^3(X,\Gm)\to H^3(C,\Gm)$ et la
fl\`eche $H^3(C,\Gm)\to H^2(k,\Pic(\Cbarre))$ issue de la suite spectrale de
Hochschild--Serre (compte tenu que $H^q(\Cbarre,\Gm)=0$ pour $q\geq 2$).
\end{lem}

\begin{proof}
Quitte \`a remplacer~$X$ par~$C$ et~$D$ par $C \cap D$, on peut supposer que~$X$
est une courbe et que $C=X$. Notons $i : D \hookrightarrow C$ et
$j : U \hookrightarrow C$ les injections canoniques.
Par construction, l'application~$\theta$ est la compos\'ee de l'isomorphisme
$H^1(D,\Q/\Z) \isoto H^2(D,\Z)$ et de l'application $H^2(D,\Z) \to H^3(C,\Gm)$
bord de la suite exacte de faisceaux \'etales sur~$C$
\begin{align}
\label{se:divweil}
\xymatrix{
0 \ar[r] & \Gm \ar[r] & j_* \Gm \ar[r] & i_* \Z \ar[r] & 0 \rlap{\text{.}}
}
\end{align}
Notant $b : i_*\Z \to \Gm[1]$ le morphisme d\'efini par~(\ref{se:divweil}) dans
la cat\'egorie d\'eriv\'ee des faisceaux \'etales en groupes ab\'eliens sur~$C$,
le morphisme~$\theta$ est donc obtenu, modulo l'identification
$H^2(D,\Z) \simeq H^1(D,\Q/\Z)$, en appliquant successivement les foncteurs
$\RGamma(\Cbarre,-)$ et $H^2(k,-)$ \`a~$b$.

Comme $\RGamma(\Cbarre,i_*\Z)=\Z^{\Dbarre}$, la compos\'ee de
$\RGamma(\Cbarre,b) : \RGamma(\Cbarre,i_*\Z) \to \RGamma(\Cbarre,\Gm)[1]$, du
morphisme de troncation
$$\RGamma(\Cbarre,\Gm)[1] \to \left(\tau_{\geq 1}\RGamma(\Cbarre,\Gm)\right)[1]=\Pic(\Cbarre)$$
et du degr\'e $\Pic(\Cbarre) \to \Z$ est une fl\`eche $t : \Z^{\Dbarre} \to \Z$
entre complexes concentr\'es en degr\'e~$0$.  Elle provient donc d'une unique
application entre modules galoisiens.  Celle-ci n'est autre que l'application
induite par~$t$ entres les objets de cohomologie de degr\'e~$0$ des complexes
consid\'er\'es.  Autrement dit~$t$ est \'egale \`a la compos\'ee de l'application
$\Gamma(\Cbarre,b) : \Z^{\Dbarre} \to \Pic(\Cbarre)$ bord
de~(\ref{se:divweil}) et du degr\'e $\Pic(\Cbarre)\to\Z$. Par cons\'equent~$t$ est
l'application $m \mapsto \sum_{P\in \Dbarre}m(P)$ et
$H^2(k,t) : H^2(D,\Z)\to H^2(k,\Z)$
co\"incide donc avec~$\sigma_C$ modulo les isomorphismes canoniques
$H^2(D,\Z)\simeq H^1(D,\Q/\Z)$ et $H^2(k,\Z)\simeq H^1(k,\Q/\Z)$.  Ceci ach\`eve
la d\'emonstration puisque $H^2(k,t)$ est la compos\'ee de~$\theta$ avec les
fl\`eches indiqu\'ees dans l'\'enonc\'e du lemme.
\end{proof}

L'int\'er\^et du lemme~\ref{lemme:controleresidufort} ne se r\'esume pas au cas
o\`u~$C$ et~$D$ se rencontrent en un unique point. Notamment, une cons\'equence
imm\'ediate de ce lemme est la description du groupe de Brauer du compl\'ementaire
d'une courbe plane lisse sur un corps de caract\'eristique~$0$ arbitraire:

\begin{prop}
\label{prop:description-br}
Soit $D \subset \P^2_k$ une courbe plane lisse sur un corps~$k$ de
caract\'eristique~$0$.  Notons $U=\P^2_k\setminus D$.  Le groupe $\Br(U)$
s'inscrit dans une suite exacte canonique
\begin{align}
\label{se:loc-U}
\xymatrix{
0 \ar[r]& \Br(k) \ar[r] & \Br(U) \ar[r] & H^1(D,\Q/\Z) \ar[r]^\sigma & H^1(k,\Q/\Z)\text{,}
}
\end{align}
o\`u l'application~$\sigma$ v\'erifie
 $\sigma(m)=\sum_{P \in L \cap D} \Cores_{k(P)/k}(m(P))$
pour tout~$m$ et toute droite $L\subset \P^2_k$ rencontrant~$D$
transversalement en chaque point d'intersection.  (En particulier le membre de
droite de cette \'egalit\'e ne d\'epend-il pas du choix de~$L$.)
\end{prop}

\begin{proof}
Soit~$\kbar$ une cl\^oture alg\'ebrique de~$k$. Comme $H^q(\P^2_{\kbar},\Gm)$
s'annule pour $q \geq 2$ et comme $H^2(k,\Pic(\P^2_{\kbar})) \simeq
H^1(k,\Q/\Z)$, la suite spectrale de Hochschild--Serre induit une suite
exacte
\begin{align}
\label{se:h3p2}
\xymatrix{
H^3(k,\Gm) \ar[r] & H^3(\P^2_k,\Gm) \ar[r] & H^1(k,\Q/\Z) \rlap{\text{.}}
}
\end{align}
Par ailleurs, la suite~(\ref{se:loc}) se prolonge \`a droite en une suite exacte
\begin{align}
\label{se:loc-droite}
\xymatrix{
H^1(D,\Q/\Z) \ar[r]^{\theta} & H^3(\P^2_k,\Gm) \ar[r] & H^3(U,\Gm) \rlap{\text{.}}
}
\end{align}
La compos\'ee de la fl\`eche de gauche de~(\ref{se:h3p2}) et de la fl\`eche de
droite de~(\ref{se:loc-droite}) est injective (tout point rationnel de~$U$
en fournit une r\'etraction).  Par cons\'equent la fl\`eche de droite
de~(\ref{se:h3p2}) est injective sur l'image de~$\theta$.  La
suite~(\ref{se:loc}) reste donc exacte si l'on remplace son dernier terme
par $H^1(k,\Q/\Z)$; gr\^ace au lemme~\ref{lemme:controleresidufort}, la
proposition~\ref{prop:description-br} s'ensuit.
\end{proof}

Lorsque~$D$ est une conique, cette description se simplifie encore: on obtient
une suite exacte
\begin{equation}
\label{se:conique}
\xymatrix{
0 \ar[r] & \Br(k) \ar[r] & \Br(U) \ar[r] & k^*/k^{*2} \ar[r] & 0\rlap{\text{.}}
}
\end{equation}
Si $f \in k[x,y,z]$ est une forme quadratique s'annulant sur~$D$ et si
$\ell \in k[x,y,z]$ d\'esigne une forme lin\'eaire arbitraire, l'alg\`ebre de
quaternions $(f/\ell^2,a)$ est non ramifi\'ee sur~$U$ et d\'efinit un rel\`evement
dans $\Br(U)$ de la classe de~$a$ dans $k^*/k^{*2}$.

\bigskip
Il serait souhaitable de disposer d'une description aussi compl\`ete du groupe
de Brauer de~$U$ dans la situation o\`u~$U$ est le compl\'ementaire d'une section
hyperplane lisse dans une surface cubique lisse que dans la situation o\`u~$U$
est le compl\'ementaire d'une courbe lisse dans le plan projectif.
Malheureusement le lemme~\ref{lemme:controleresidu} ne fournit qu'une
majoration du groupe $\Ker(\theta)$.  Il s'av\`ere n\'eanmoins que la
proposition~\ref{prop:description-br} permet dans certains cas de pallier
cette difficult\'e, comme nous allons tout de suite le constater sur l'exemple
des rev\^etements cycliques de degr\'e~$n$ du plan projectif ramifi\'es le long
d'une courbe lisse de degr\'e~$n$.

Supposons, pour simplifier, que $H^3(k,\Gm)=0$ (cette hypoth\`ese est notamment
satisfaite si~$k$ est un corps de nombres).
Soient $D \subset \P^2_k$ une courbe plane lisse, de degr\'e $n \geq 1$, et
$f \in k[x,y,z]$ un polyn\^ome homog\`ene de degr\'e~$n$ s'annulant sur~$D$.  Soit
$X \subset \P^3_k$ la surface lisse d'\'equation $f(x,y,z)=t^n$.  Voyant
indiff\'eremment la courbe~$D$ comme plong\'ee dans~$\P^2_k$ ou dans~$X$ (auquel
cas elle s'identifie \`a la section hyperplane d'\'equation $t=0$), nous posons
$U=X\setminus D$ et $V=\P^2_k \setminus D$.  Le rev\^etement $\pi:X \to \P^2_k$
d\'efini par $[x:y:z:t]\mapsto [x:y:z]$ est ramifi\'e \`a l'ordre~$n$ le long de~$D$
et induit sur~$U$ une structure de torseur sous~$\mmu_n$ au-dessus de~$V$.
Par cons\'equent les suites exactes~(\ref{se:loc}) et~(\ref{se:loc-U}) associ\'ees
aux immersions ouvertes $U \subset X$ et $V \subset \P^2_k$ s'inscrivent dans
un diagramme commutatif
\begin{align}
\begin{aligned}
\label{diag:grand}
\xymatrix{
0 \ar[r] & \Br(k) \ar[r]\ar[d] & \Br(V) \ar[r]\ar[d]^{\pi^*} & H^1(D,\Q/\Z) \ar[r]^\sigma\ar[d]^{\times n} & H^1(k,\Q/\Z)\ar[d]\\
0 \ar[r] & \Br(X) \ar[r]\ar[d] & \Br(U) \ar[r]\ar[d]^{\pi_*} & H^1(D,\Q/\Z) \ar[r]^\theta\ar@{=}[d] & H^3(X,\Gm)\ar[d]\\
0 \ar[r] & \Br(k) \ar[r] & \Br(V) \ar[r] & H^1(D,\Q/\Z) \ar[r]^\sigma & H^1(k,\Q/\Z)
}
\end{aligned}
\end{align}
(cf.~\cite[Prop.~1.1.1, Prop.~1.1.2]{CTSD94}; remarquer que
$H^3(\P^2_k,\Gm)=H^1(k,\Q/\Z)$ puisque $H^3(k,\Gm)=0$).
Il r\'esulte de la partie
inf\'erieure de ce diagramme que $\Ker(\theta) \subset \Ker(\sigma)$, tandis que
la partie sup\'erieure montre que tout \'el\'ement de $\Ker(\sigma)$ d'ordre premier
\`a~$n$ appartient \`a $\Ker(\theta)$.  Ainsi avons-nous \'etabli~le

\begin{lem}
Avec les notations ci-dessus, pour tout $m \in H^1(D,\Q/\Z)$ d'ordre premier
\`a~$n$, on a $m \in \Ker(\theta)$ si et seulement si $m \in \Ker(\sigma)$.
\end{lem}

Il s'agit l\`a d'un crit\`ere simple pour l'appartenance de~$m$ \`a $\Ker(\theta)$,
lorsque l'ordre de~$m$ est premier \`a~$n$.  Quand ce crit\`ere s'applique, le
diagramme~(\ref{diag:grand}) montre de plus que pour exhiber une classe
de~$\Br(U)$ dont le r\'esidu le long de~$D$ soit \'egal \`a~$m$, il suffit de
chercher une telle classe dans l'image de~$\pi^*$.  C'est en proc\'edant ainsi
que nous avons explicit\'e l'alg\`ebre de quaternions~$B'$ de la
proposition~\ref{prop:brauer-affine2} \`a partir de la classe dans
$H^1(D,\Q/\Z)$ de la $2$\nobreakdash-isog\'enie~(\ref{se:m2}) (dont il n'\'etait
pas clair \emph{a priori} qu'elle appartenait au noyau de~$\theta$),
\`a l'aide de la proposition~\ref{prop:description-br}.

\medskip
Mentionnons pour terminer ce paragraphe une derni\`ere application du
lemme~\ref{lemme:controleresidufort}:

\begin{prop}
\label{prop:finitude-Br-U}
Soient~$X$ une vari\'et\'e propre et lisse sur un corps de nombres~$k$ et
$D \subset X$ une sous-vari\'et\'e ferm\'ee de codimension~$1$, lisse et
g\'eom\'etriquement connexe sur~$k$.
Notons~$\kbar$ une cl\^oture alg\'ebrique de~$k$ et posons $U=X\setminus D$.  S'il
existe une courbe propre et lisse $C \subset X \otimes_k\kbar$ rencontrant
$D \otimes_k\kbar$ en un unique point, avec intersection transverse en ce
point, le quotient $\Br(U)/\Br(X)$ est fini.
\end{prop}

\begin{proof}
Nous devons montrer que le noyau de l'application~$\theta$ apparaissant
dans~(\ref{se:loc}) est fini.  Pour ceci il est loisible de remplacer~$k$ par
une extension finie arbitraire puisque le noyau de la fl\`eche de restriction
$H^1(D,\Q/\Z) \to H^1(D \otimes_k \ell,\Q/\Z)$ est fini pour toute extension
finie $\ell/k$.  En particulier pouvons-nous supposer la courbe~$C$ d\'efinie
sur~$k$.  D'apr\`es le lemme~\ref{lemme:controleresidu} il suffit alors de
s'assurer que $\Ker(\sigma_C)$ est fini.  Or $\sigma_C$ est ici une r\'etraction
de la fl\`eche naturelle $H^1(k,\Q/\Z) \to H^1(D,\Q/\Z)$.  Le noyau
de~$\sigma_C$ s'identifie donc au groupe $H^0(k,H^1(\Dbarre,\Q/\Z))$, lequel est
fini d'apr\`es Katz et Lang (voir la remarque~\ref{rq:HSC}~(i)).
\end{proof}

Les hypoth\`eses de la proposition~\ref{prop:finitude-Br-U} sont notamment
satisfaites lorsque~$X$ est une sous-vari\'et\'e de~$\P^n_k$ de dimension~$\geq 2$
contenant g\'eom\'etriquement une droite et que~$D$ est une section hyperplane
lisse de~$X$ ne contenant pas cette droite (par exemple~$X$ pourrait \^etre une
surface cubique lisse et~$D$ une section hyperplane lisse arbitraire).

Que l'on ne puisse se dispenser de supposer l'existence de~$C$ se voit d\'ej\`a
sur l'exemple du compl\'ementaire d'une conique dans le plan, compte tenu de la
suite exacte~(\ref{se:conique}).

\subsection{Contre-exemples \`a l'approximation forte}

Il r\'esulte du lemme~\ref{lemme:points-locaux} et de la
proposition~\ref{prop:evalBpdyadique} que quel que soit l'entier $a \neq 0$,
la classe~$B'$ de la proposition~\ref{prop:brauer-affine2} est \emph{toujours}
responsable d'une obstruction de Brauer--Manin \`a l'approximation forte sur la
surface affine~$U'$ sur~$\Q$ d'\'equation $x^3+y^3+2z^3=a$ (et m\^eme \`a
l'approximation forte en dehors de la place r\'eelle, puisque~$B'$ s'\'evalue
trivialement sur $U'(\R)$; voir la remarque~\ref{rq:evalB}).

\begin{exemple}
\label{ex:approx-forte}
L'\'equation $x^3+y^3+2z^3=2$ n'admet pas de solution $x,y,z \in \Z$ telle que
$\congru{x+y}{2}{8}$ et $\congru{z}{2}{4}$, bien qu'une solution dans~$\Z_2$
satisfaisant ces congruences existe (par exemple $x_0=9$, $y_0=9$ et
$z_0=-2\sqrt[3]{91}$) et bien qu'il existe des solutions dans~$\Q$
arbitrairement proches de toute solution dyadique fix\'ee (par exemple la
solution $x=-1/15$, $y=-17/15$, $z=6/5$ v\'erifie $v_2(x-x_0) \geq 3$,
$v_2(y-y_0) \geq 3$ et $v_2(z-z_0)\geq 2$).
\end{exemple}

\begin{proof}
La surface cubique projective~$X'$ d'\'equation $x^3+y^3+2z^3=2t^3$ est
rationnelle sur~$\Q$ puisqu'elle contient deux droites gauches conjugu\'ees
(cf.~\cite{swdbirationality}).
Par cons\'equent elle satisfait l'approximation
faible; en particulier l'ensemble $U'(\Q)$ est dense dans $U'(\Q_2)$.
Supposons maintenant qu'il existe un point $M \in \sU'(\Z)$ de coordonn\'ees
$(x,y,z)$ tel que $\congru{x+y}{2}{8}$ et $\congru{z}{2}{4}$.
Notant $B' \in \Br(U')$ la classe d\'efinie dans l'\'enonc\'e de la
proposition~\ref{prop:brauer-affine2} (avec $a=2$), l'invariant dyadique de
$B'(M)$ est alors non nul (cf.~\cite[Ch.~III, th.~1]{serrecoursarithmetique}).
Or, d'apr\`es la remarque~\ref{rq:evalB}, l'\'evaluation de~$B'$ sur $\sU'(\Z_p)$
est triviale pour tout $p \neq 2$.  La loi de r\'eciprocit\'e globale appliqu\'ee \`a
$B'(M) \in \Br(\Q)$ fournit donc une contradiction.
\end{proof}

\begin{rmq}
\label{rq:cassels}
D'apr\`es Cassels~\cite{casselsmathcomp}, si $x,y,z \in \Z$ v\'erifient $x^3+y^3+z^3=3$,
alors $\congru{x \equiv y}{z}{9}$.
Notons comme pr\'ec\'edemment~$U$ la surface affine d'\'equation $x^3+y^3+z^3=3$
et~$X$ la surface cubique projective correspondante.
Les propositions~\ref{prop:proj-br} et~\ref{prop:brauer-affine}
montrent que le groupe $\Br(U)/\Br(\Q)$ est d'ordre~$3$, engendr\'e par l'image
de la classe $A \in \Br(X)$ d\'efinie dans l'\'enonc\'e de la
proposition~\ref{prop:proj-br}.  Comme~$X$ et~$A$ ont bonne r\'eduction hors
de~$3$, l'\'evaluation de~$A$ sur $X(\Q_p)$ est nulle pour tout nombre premier
$p\neq 3$.  L'\'evaluation de~$A$ sur $X(\R)$ \'etant \'egalement nulle, il r\'esulte
de la loi de r\'eciprocit\'e globale que l'invariant $3$\nobreakdash-adique de
$A(M)$ s'annule pour tout $M \in X(\Q)$.  Concr\`etement, cela signifie que si
$x,y,z \in \Q$ v\'erifient $x^3+y^3+z^3=3$, alors l'invariant local de
$(x+jy,3)_j \in \Br(\Q(j))$ en l'unique place de $\Q(j)$ divisant~$3$ est nul;
on constate \`a l'aide du formulaire~\cite[p.~34]{CTKS} que l'affirmation de
Cassels en r\'esulte.  Cette interpr\'etation du calcul de Cassels montre bien que
le d\'efaut d'approximation forte mis en \'evidence dans~\cite{casselsmathcomp}
est en r\'ealit\'e d\^u \`a une obstruction de Brauer--Manin \`a l'approximation faible:
l'ensemble $U(\Q) \cap \sU(\Z_3)$ n'est m\^eme pas dense dans $\sU(\Z_3)$.
L'exemple~\ref{ex:approx-forte}, en revanche, pr\'esente indiscutablement un
d\'efaut d'approximation forte puisque l'approximation faible y est satisfaite.
\end{rmq}

\subsection{Crit\`eres pour l'existence de points entiers}

Le th\'eor\`eme~\ref{thm:sur-les-entiers} am\`ene naturellement les deux questions
suivantes:

\begin{questions}
\label{qu}
(i) Soit $f \in \Z[x,y,z]$ un polyn\^ome homog\`ene de degr\'e~$3$ tel que la courbe
plane d'\'equation $f=0$ soit lisse sur~$\Q$.  Soit~$n$ un entier non nul.
L'\'equation $f(x,y,z)=n$ admet-elle une solution $(x,y,z)\in \Z^3$ d\`es que
l'obstruction de Brauer--Manin enti\`ere ne s'y oppose pas~?

(ii) Plus g\'en\'eralement, soient $X \subset \P^3_k$ une surface cubique lisse sur
un corps de nombres~$k$ et $U \subset X$ le compl\'ementaire d'une section
hyperplane lisse.  L'ensemble $U(k)$ est-il dense dans
$U(\A_k)_\bullet^{\Br(U)}$~?
\end{questions}

La notation $U(\A_k)_\bullet^{\Br(U)}$ d\'esigne ici le sous-ensemble de
$U(\A_k^f) \times \pi_0(U(\A_k^\infty))$ constitu\'e des familles orthogonales \`a
$\Br(U)$ pour l'accouplement de Brauer--Manin, o\`u $\A_k=\A_k^f \times
\A_k^\infty$ est la d\'ecomposition des ad\`eles de~$k$ en ad\`eles finis et
infinis.

Si~$U$ est le compl\'ementaire d'une section hyperplane lisse dans une surface
cubique lisse sur un corps de nombres, une r\'eponse affirmative \`a la
question~(ii) aurait pour corollaire que l'ensemble des points entiers de~$U$
(un mod\`ele \'etant fix\'e) serait soit vide, soit dense dans~$U$ pour la topologie
de Zariski (en particulier infini), compte tenu de la
proposition~\ref{prop:finitude-Br-U}.  Beukers~\cite{beukers} (voir aussi
\cite[Th.~6.13]{hassetttschinkel}) a \'etabli la Zariski-densit\'e potentielle des
points entiers de~$U$ (c'est-\`a-dire la Zariski-densit\'e apr\`es une extension
finie des scalaires).  Hassett et
Tschinkel~\cite[\textsection7]{hassetttschinkel} discutent plus g\'en\'eralement
la question de la Zariski-densit\'e potentielle des points entiers sur les
surfaces log\nobreakdash-K3.  Il~serait cependant d\'eraisonnable de poser la
question~\ref{qu}~(ii) pour toutes les surfaces log\nobreakdash-K3. La r\'eponse
est en effet d\'ej\`a n\'egative pour les surfaces log-del Pezzo, comme l'illustrent
les deux exemples ci-dessous.

\begin{exemple}
\label{exemple1}
De fa\c{c}on \'evidente, l'\'equation $2x^2+3y^2+4z^2=1$ n'admet pas de solution enti\`ere
$(x,y,z)\in\Z^3$.  Nous allons voir qu'il n'y a pourtant pas d'obstruction de
Brauer--Manin \`a l'existence de telles solutions.

Notons $X \subset \P^3_\Q$ la surface quadrique d\'efinie par
$2x^2+3y^2+4z^2=t^2$ et $D \subset X$ la section hyperplane lisse d'\'equation
$t=0$.  Soient de plus $\sU=\Spec(\Z[x,y,z]/(2x^2+3y^2+4z^2-1))$ et
$U=X \setminus D=\sU\otimes_\Z\Q$.  Le lemme~\ref{lemme:controleresidufort} et
la suite exacte~(\ref{se:loc}) permettent de v\'erifier que $\Br(U)/\Br(\Q)$ est
d'ordre~$2$, engendr\'e par la classe de l'alg\`ebre de quaternions $A=(1-2z,-6)$.
Notons $P,Q\in U(\Q)$ les points de coordonn\'ees $P=(0,1/2,1/4)$ et
$Q=(3/5,1/5,1/5)$.  Vu les d\'enominateurs des coordonn\'ees de~$P$ et de~$Q$, on
a $\sU(\Z_p)\neq\emptyset$ pour tout~$p$.  D'autre part~$P$ et~$Q$ sont tous
les deux entiers en la place~$3$ et les invariants $3$\nobreakdash-adiques
de~$A(P)$ et de~$A(Q)$ sont diff\'erents.  Il existe donc bien une famille de
points entiers locaux de~$\sU$ orthogonale \`a $\Br(U)$.
\end{exemple}

Le premier auteur et Xu~\cite[Th.~6.3]{ctxu} d\'emontrent n\'eanmoins que si~$n$ est un entier non nul
et $f \in \Z[x,y,z]$ une forme quadratique
homog\`ene \emph{ind\'efinie}, alors l'\'equation $f(x,y,z)=n$ admet une solution
enti\`ere d\`es que l'obstruction de Brauer--Manin ne s'y oppose pas.
Pour autant, m\^eme avec~$f$ ind\'efinie, la surface log\nobreakdash-del Pezzo
compl\'ementaire dans le plan projectif de la conique lisse d'\'equation $f=0$
peut ne pas admettre de point entier sans que cela r\'esulte d'une
obstruction de Brauer--Manin enti\`ere:

\begin{exemple}
\label{exemple2}
Soit $f(x,y,z)=16x^2 + 9y^2 - 3z^2$. Notons $\sD \subset \P^2_\Z$ le ferm\'e
d\'efini par $f=0$, puis $\sV=\P^2_\Z \setminus \sD$ et $V=\sV\otimes_\Z\Q$.
Notons de plus $X \subset \P^3_\Q$ la surface quadrique d'\'equation
$f(x,y,z)=t^2$. Posons $\sU=\Spec(\Z[x,y,z]/(f(x,y,z)-1))$, de sorte que
$U=\sU\otimes_\Z\Q$ s'identifie au compl\'ementaire dans~$X$ de la section
hyperplane lisse d'\'equation $t=0$ et est naturellement muni d'un morphisme
fini \'etale $\pi:U\to V$ de degr\'e~$2$.  L'\'equation $f(x,y,z)=-1$ n'admet pas de
solution $(x,y,z)\in \Z^3$ car elle n'en admet pas dans $(\Z/3\Z)^3$.  Par
ailleurs, d'apr\`es \cite[Ex.~1.2]{schulzepillotxu} ou \cite[Prop.~8.2]{ctxu},
l'\'equation $f(x,y,z)=1$ n'admet pas de solution $(x,y,z)\in\Z^3$ mais en admet
une dans~$\R^3$ et dans~$\Z_p^3$ pour tout~$p$ (autrement dit $\sU$ n'admet
pas de point sur~$\Z$ mais admet des points entiers locaux).  Comme~$1$
et~$-1$ sont les seules unit\'es de~$\Z$, on conclut que $\sV(\Z)=\emptyset$.
D'autre part, la suite exacte~(\ref{se:conique}) montre que le groupe
$\Br(V)/\Br(\Q)$ est annul\'e par~$2$ et le diagramme~(\ref{diag:grand}) permet
d'en d\'eduire que $\pi^*\Br(V) \subset \Br(X)$; comme $\Br(X)=\Br(\Q)$, il
s'ensuit que la projection dans~$\sV$ de n'importe quelle famille de points
locaux de~$\sU$ est orthogonale \`a $\Br(V)$ pour l'accouplement de
Brauer--Manin.  Il n'y a donc pas d'obstruction de Brauer--Manin enti\`ere \`a
l'existence d'une solution dans~$\Z^3$ de l'\'equation
$16x^2 + 9y^2 - 3z^2=\pm 1$, bien qu'une telle solution n'existe pas.

Le lecteur v\'erifiera que cet exemple se g\'en\'eralise \`a
$f(x,y,z)=16m^2x^2+p^{2k}y^2-pz^2$ avec $m$, $k>0$ entiers et~$p$ premier
congru \`a~$3$ modulo~$8$.
\end{exemple}

L'absence de point entier sur le sch\'ema~$\sU$ de l'exemple~\ref{exemple1} est
li\'ee \`a un ph\'enom\`ene archim\'edien.  Dans l'exemple~\ref{exemple2}, bien qu'il
n'y ait pas d'obstruction de Brauer--Manin \`a l'existence d'un point
entier sur~$\sU$, la vacuit\'e de $\sU(\Z)$ s'explique par une obstruction de
Brauer--Manin enti\`ere sur un rev\^etement \'etale de~$U$.

Ces deux difficult\'es ne concernent vraisemblablement pas la question~\ref{qu}.
En effet celle-ci porte sur une \'equation de degr\'e impair, ce qui \'ecarte tout
ph\'enom\`ene archim\'edien \'evident, et d'autre part, d'apr\`es~\cite{serreabhyankar},
la surface~$U$ qui y est consid\'er\'ee est simplement connexe.

Examinons, pour terminer, la question~\ref{qu}~(ii) dans la situation de
l'exemple~\ref{ex:approx-forte}.

\begin{exemple}
Posons $\sU' = \Spec(\Z[x,y,z]/(x^3+y^3+2z^3-2))$.  Nous avons vu que la
surface $U'=\sU' \otimes_\Z \Q$ est rationnelle sur~$\Q$; la
proposition~\ref{prop:brauer-affine2} entra\^ine donc que $\Br(U')/\Br(\Q)$ est
d'ordre~$2$, engendr\'e par l'image de~$B'$.  D'autre part, l'\'evaluation de~$B'$
sur $U'(\R)$ et sur $\sU'(\Z_p)$ pour tout $p \neq 2$ est triviale
(remarque~\ref{rq:evalB}) et~$B'$ atteint la valeur~$0$ sur $\sU'(\Z_2)$
(proposition~\ref{prop:evalBpdyadique}).  Pour que la question~\ref{qu}~(ii)
admette une r\'eponse positive, il est donc n\'ecessaire que
$\sU'(\Z)$ soit dense dans $\prod_{p \notin \{2,\infty\}} \sU'(\Z_p)$
pour la topologie produit.
\end{exemple}

\bibliographystyle{smfalpha}
\bibliography{troiscubes}

\providecommand{\bysame}{\leavevmode ---\ }
\providecommand{\og}{``}
\providecommand{\fg}{''}
\providecommand{\smfandname}{et}
\providecommand{\smfedsname}{\'eds.}
\providecommand{\smfedname}{\'ed.}
\providecommand{\smfmastersthesisname}{M\'emoire}
\providecommand{\smfphdthesisname}{Th\`ese}
\begin{thebibliography}{CTKS87}

\bibitem[Beu99]{beukers}
{\scshape F.~Beukers} -- {\og Integral points on cubic surfaces\fg}, Number
  theory ({O}ttawa, {ON}, 1996), CRM Proc. Lecture Notes, vol.~19, Amer. Math.
  Soc., Providence, RI, 1999, p.~25--33.

\bibitem[Cas85]{casselsmathcomp}
{\scshape J.~W.~S. Cassels} -- {\og A note on the {D}iophantine equation {$x\sp
  3+y\sp 3+z\sp 3=3$}\fg}, \emph{Math. Comp.} \textbf{44} (1985), no.~169,
  p.~265--266.

\bibitem[Cre97]{cremona}
{\scshape J.~E. Cremona} -- \emph{Algorithms for modular elliptic curves},
  seconde \smfedname, Cambridge University Press, Cambridge, 1997.

\bibitem[CTKS87]{CTKS}
{\scshape J.-L. Colliot-Th{\'e}l{\`e}ne, D.~Kanevsky {\normalfont \smfandname}
  J.-J. Sansuc} -- {\og Arithm\'etique des surfaces cubiques diagonales\fg},
  Diophantine approximation and transcendence theory ({B}onn, 1985), Lecture
  Notes in Math., vol. 1290, Springer, Berlin, 1987, p.~1--108.

\bibitem[CTS87]{CTS}
{\scshape J.-L. Colliot-Th{\'e}l{\`e}ne {\normalfont \smfandname} J.-J. Sansuc}
  -- {\og La descente sur les vari\'et\'es rationnelles, {II}\fg}, \emph{Duke
  Math. J.} \textbf{54} (1987), no.~2, p.~375--492.

\bibitem[CTSD94]{CTSD94}
{\scshape J.-L. Colliot-Th{\'e}l{\`e}ne {\normalfont \smfandname}
  P.~Swinnerton-Dyer} -- {\og Hasse principle and weak approximation for
  pencils of {S}everi-{B}rauer and similar varieties\fg}, \emph{J. reine angew.
  Math.} \textbf{453} (1994), p.~49--112.

\bibitem[CTX09]{ctxu}
{\scshape J.-L. Colliot-Th{\'e}l{\`e}ne {\normalfont \smfandname} F.~Xu} --
  {\og Brauer--{M}anin obstruction for integral points of homogeneous spaces
  and representation by integral quadratic forms\fg}, \emph{Compositio Math.}
  \textbf{145} (2009), p.~309--363.

\bibitem[CV94]{connvaserstein}
{\scshape W.~Conn {\normalfont \smfandname} L.~N. Vaserstein} -- {\og On sums
  of three integral cubes\fg}, The {R}ademacher legacy to mathematics
  ({U}niversity {P}ark, {PA}, 1992), Contemp. Math., vol. 166, Amer. Math.
  Soc., Providence, RI, 1994, p.~285--294.

\bibitem[EJ09]{elsenhansjahnel}
{\scshape A.-S. Elsenhans {\normalfont \smfandname} J.~Jahnel} -- {\og New sums
  of three cubes\fg}, \emph{Math. Comp.} \textbf{78} (2009), no.~266,
  p.~1227--1230.

\bibitem[Fuj02]{fujiwara}
{\scshape K.~Fujiwara} -- {\og A proof of the absolute purity conjecture (after
  {G}abber)\fg}, Algebraic geometry 2000, {A}zumino ({H}otaka), Adv. Stud. Pure
  Math., vol.~36, Math. Soc. Japan, Tokyo, 2002, p.~153--183.

\bibitem[Gro68]{brauerIII}
{\scshape A.~Grothendieck} -- {\og Le groupe de {B}rauer {III} : exemples et
  compl\'ements\fg}, Dix expos\'es sur la cohomologie des sch\'emas, Advanced
  studies in pure mathematics, vol.~3, Masson \& North-Holland, Paris,
  Amsterdam, 1968, p.~88--188.

\bibitem[GS06]{gilleszamuely}
{\scshape P.~Gille {\normalfont \smfandname} T.~Szamuely} -- \emph{Central
  simple algebras and {G}alois cohomology}, Cambridge Studies in Advanced
  Mathematics, vol. 101, Cambridge University Press, Cambridge, 2006.

\bibitem[Guy04]{guy}
{\scshape R.~K. Guy} -- \emph{Unsolved problems in number theory}, troisi\`eme
  \smfedname, Problem Books in Mathematics, Springer-Verlag, New York, 2004.

\bibitem[HB92]{HBdensity}
{\scshape D.~R. Heath-Brown} -- {\og The density of zeros of forms for which
  weak approximation fails\fg}, \emph{Math. Comp.} \textbf{59} (1992), no.~200,
  p.~613--623.

\bibitem[HT01]{hassetttschinkel}
{\scshape B.~Hassett {\normalfont \smfandname} {\relax Yu}.~Tschinkel} -- {\og
  Density of integral points on algebraic varieties\fg}, Rational points on
  algebraic varieties, Progr. Math., vol. 199, Birkh\"auser, Basel, 2001,
  p.~169--197.

\bibitem[HW08]{hardywright}
{\scshape G.~H. Hardy {\normalfont \smfandname} E.~M. Wright} -- \emph{An
  introduction to the theory of numbers}, sixi\`eme \smfedname, Oxford
  University Press, Oxford, 2008.

\bibitem[KL81]{katzlang}
{\scshape N.~M. Katz {\normalfont \smfandname} S.~Lang} -- {\og Finiteness
  theorems in geometric class field theory\fg}, \emph{L'Enseign. Math. (2)}
  \textbf{27} (1981), no.~3-4, p.~285--319.

\bibitem[Ko36]{ko}
{\scshape C.~Ko} -- {\og Decompositions into four cubes\fg}, \emph{J. Lond.
  Math. Soc.} \textbf{11} (1936), p.~218--219.

\bibitem[Koy00]{koyama}
{\scshape K.~Koyama} -- {\og On searching for solutions of the {D}iophantine
  equation {$x\sp 3+y\sp 3+2z\sp 3=n$}\fg}, \emph{Math. Comp.} \textbf{69}
  (2000), no.~232, p.~1735--1742.

\bibitem[KT08]{KT}
{\scshape A.~Kresch {\normalfont \smfandname} {\relax Yu}.~Tschinkel} -- {\og
  Two examples of {B}rauer--{M}anin obstruction to integral points\fg},
  \emph{Bull. Lond. Math. Soc.} \textbf{40} (2008), no.~6, p.~995--1001.

\bibitem[Man72]{Manin}
{\scshape {\relax Yu}.~I. Manin} -- \emph{Cubic forms: algebra, geometry,
  arithmetic}, Nauka, Moscou, 1972 (en russe), trad. anglaise: Cubic forms,
  North-Holland, Amsterdam, 1974, seconde \'ed., 1986.

\bibitem[Mil80]{milne}
{\scshape J.~S. Milne} -- \emph{\'{E}tale cohomology}, Princeton Mathematical
  Series, vol.~33, Princeton University Press, Princeton, N.J., 1980.

\bibitem[Ryl25]{ryley}
{\scshape S.~Ryley} -- \emph{The Ladies' Diary} (1825), p.~35, question~1420.

\bibitem[SD70]{swdbirationality}
{\scshape H.~P.~F. Swinnerton-Dyer} -- {\og The birationality of cubic surfaces
  over a given field\fg}, \emph{Michigan Math. J.} \textbf{17} (1970),
  p.~289--295.

\bibitem[Sel51]{Selmer}
{\scshape E.~S. Selmer} -- {\og The {D}iophantine equation {$ax\sp 3+by\sp
  3+cz\sp 3=0$}\fg}, \emph{Acta Math.} \textbf{85} (1951), p.~203--362.

\bibitem[Ser62]{corpslocaux}
{\scshape J.-P. Serre} -- \emph{Corps locaux}, Publications de l'Institut de
  Math\'ematique de l'Universit\'e de Nancago, VIII, Actualit\'es Sci. Indust.,
  No. 1296. Hermann, Paris, 1962.

\bibitem[Ser70]{serrecoursarithmetique}
\bysame , \emph{Cours d'arithm\'etique}, collection SUP: \og{}le
  math\'ematicien\fg{}, vol.~2, Presses Universitaires de France, Paris, 1970.

\bibitem[Ser94]{cohomologiegaloisienne}
\bysame , \emph{Cohomologie galoisienne}, Cinqui\`eme \'edition, Lecture Notes
  in Mathematics, vol.~5, Springer-Verlag, Berlin, 1994.

\bibitem[Ser95]{serreabhyankar}
\bysame , {\og Rev\^etements ramifi\'es du plan projectif (d'apr\`es
  {S}.~{A}bhyankar)\fg}, S\'eminaire {B}ourbaki, vol.\ 5, Soc. Math. France,
  Paris, 1995, p.~483--489, exp.~204.

\bibitem[SPX04]{schulzepillotxu}
{\scshape R.~Schulze-Pillot {\normalfont \smfandname} F.~Xu} -- {\og
  Representations by spinor genera of ternary quadratic forms\fg}, Algebraic
  and arithmetic theory of quadratic forms, Contemp. Math., vol. 344, Amer.
  Math. Soc., Providence, RI, 2004, p.~323--337.

\bibitem[Tor38]{tornheim}
{\scshape L.~Tornheim} -- {\og Sums of {$n$}-th powers in fields of prime
  characteristic\fg}, \emph{Duke Math. J.} \textbf{4} (1938), no.~2,
  p.~359--362.

\bibitem[Vas91]{vaserstein}
{\scshape L.~N. Vaserstein} -- {\og Sums of cubes in polynomial rings\fg},
  \emph{Math. Comp.} \textbf{56} (1991), no.~193, p.~349--357.

\bibitem[V{\'e}l71]{velu}
{\scshape J.~V{\'e}lu} -- {\og Isog\'enies entre courbes elliptiques\fg},
  \emph{C. R. Acad. Sci. Paris S\'er. A-B} \textbf{273} (1971), p.~A238--A241.

\end{thebibliography}
\end{document}